\documentclass[11pt]{amsart}
\usepackage{amssymb}
\usepackage[cmtip,matrix,arrow]{xy}

\newtheorem{theorem}{Theorem}[section]
\newtheorem{corollary}[theorem]{Corollary}
\newtheorem{proposition}[theorem]{Proposition}

\newtheorem{lemma}[theorem]{Lemma}
\theoremstyle{remark}
\newtheorem{remark}[theorem]{Remark}

\newtheorem{remarks}[theorem]{Remarks}
\newtheorem{example}[theorem]{Example}

\newcommand\be{\begin{equation}\label}
\newcommand\ee{\end{equation}}
\newcommand\M{\mathcal{M}}

\renewcommand{\O}{\mathcal{O}}
\newcommand{\Co}{\mathcal{C}}

\newcommand{\U}{\on{U}}

\newcommand{\F}{\mathcal{F}}

\newcommand{\R}{\mathbb{R}}
\newcommand{\C}{\mathbb{C}}

\newcommand{\Z}{\mathbb{Z}}

\newcommand{\pr}{\on{pr}}

\newcommand\lie[1]{\mathfrak{#1}}
\renewcommand{\k}{\lie{k}}

\newcommand{\g}{\lie{g}}
\newcommand{\m}{\lie{m}}
\newcommand{\n}{\lie{n}}

\renewcommand{\t}{\lie{t}}
\newcommand{\Alc}{\lie{A}}

\newcommand{\on}{\operatorname}

\newcommand{\Ad}{ \on{Ad} } 
\newcommand{\ad}{ \on{ad} } 

\newcommand{\Eul}{ \on{Eul} }

 \newcommand{\Ind}{ \on{Ind}}
\renewcommand{\ker}{ \on{ker}} 
 
\newcommand{\SU}{ \on{SU}}

\newcommand{\Vol}{ \on{Vol}}

\newcommand{\D}{ \mathcal{D} }

\newcommand\dirac{/\kern-1.2ex\partial} 
\newcommand\qu{/\kern-.7ex/} 

\newcommand{\id}{\on{id}}

\newcommand{\hra}{\hookrightarrow}

\renewcommand{\d}{{\on{d}}}
\newcommand{\ol}{\overline}

\newcommand\Phinv{\Phi^{-1}}

\newcommand\Sig{\Sigma}
\newcommand\sig{\sigma}
\newcommand\eps{\epsilon}
\newcommand\Om{\Omega}
\newcommand\om{\omega}
\newcommand{\Del}{\Delta}

\newcommand{\f}{\frac}

\newcommand{\p}{\partial}
\renewcommand{\l}{\langle}
\renewcommand{\r}{\rangle}
\newcommand{\hh}{{\f{1}{2}}}
\newcommand{\ti}{\tilde}

\newcommand\pt{\on{pt}}

\newcommand\vol{\on{vol}}

\newcommand\beqn{\begin{equation}}      
\newcommand\eeqn{\end{equation}}      
\newcommand{\ca}{\mathcal}

\newcommand{\mf}{\mathfrak}
\newcommand{\beq}{\begin{eqnarray*}}
\newcommand{\eeq}{\end{eqnarray*}}
\newcommand{\Duf}{\on{Duf}}
\newcommand{\tpi}{2\pi i}

\begin{document}
\sloppy

\title[Intersection pairings] {Witten's formulas for
intersection pairings on moduli spaces of flat $G$-bundles} 
\date{\today}

\author{E. Meinrenken}
\address{University of Toronto, Department of Mathematics,
100 St George Street, Toronto, Ontario M5S3G3, Canada }
\email{mein@math.toronto.edu}

\begin{abstract}In a 1992 paper \cite{wi:tw}, Witten gave a formula for the 
intersection pairings of the moduli space of flat $G$-bundles over an
oriented surface, possibly with markings. In this paper, we give a
general proof of Witten's formula, for arbitrary compact, simple
groups, and any markings for which the moduli space has at most
orbifold singularities. 
\end{abstract}

\maketitle
\tableofcontents
\section{Introduction}
Let $\Sig$ be a compact, connected, oriented surface of genus $s$ with 
$r\ge 1$ boundary components, and $G$ a compact, connected Lie group. 
Given conjugacy classes 
$\Co_1,\ldots,\Co_r$ in $G$, let 
\begin{equation}\label{eq:modulispace}
 \M(\Sig;\Co_1,\ldots,\Co_r)
\end{equation}
denote the moduli space of flat $G$-bundles over $\Sig$, with holonomy
around the $j$th boundary component in the prescribed conjugacy class
$\Co_j$. For `generic' conjugacy classes $\Co_j$, the moduli space
has the structure of a smooth, compact, connected orbifold.

Of particular interest to algebraic geometry is the case $G=\SU(n)$,
$r=1$, and $\Co=\{c\}$ the conjugacy class consisting of a generator
of the center of $\SU(n)$. That is, $c=\exp(\tpi \f{d}{n})\,I$ where
$d$ and $n$ are coprime. In this case, the moduli space $\M(\Sig;c)$
is {\em smooth}, and the Narasimhan-Seshadri theorem \cite{na:st}
identifies the moduli space with a moduli space of stable holomorphic
vector bundles of rank $n$ and degree $d$ over a Riemann surface of genus 
$s$. In \cite{ha:co} Harder-Narasimhan used methods from number theory to
calculate the Poincar\'{e} polynomial of the space $\M(\Sig;c)$.  In
1984, Atiyah-Bott \cite{at:mo} gave a more geometric computation of
the Poincar\'{e} polynomial, based on the Morse theory of the Yang-Mills 
functional. Furthermore, they constructed classes
\beq
\mf{a}^p&\in& H^{2d}(\M(\Sig;c)),\\
\mf{b}^p_{i}&\in& H^{2d-1}(\M(\Sig;c)),\ \ i=1,\ldots,2h\\
\mf{f}^p&\in& H^{2d-2}(\M(\Sig;c)),
\eeq
for any invariant polynomial $p\in \on{Pol}^d(\mf{su}(n))^G$, which
generate the cohomology ring as $p$ ranges over the set of polynomials
$p(\xi)=\on{tr}(\xi^k)$, $k=2,\ldots,n$.  Beauville \cite{be:sur} gave
an alternative construction of the Atiyah-Bott classes, and
Biswas-Raghavendra \cite{bi:can} obtained generators for moduli spaces
of parabolic bundles. See also Racani\`{e}re \cite{ra:ki}.

In 1992 Witten \cite{wi:tw} proposed general formulas computing all
intersection pairings between Atiyah-Bott classes, generalizing
results of Thaddeus \cite{th:co} for the $n=2$ case. Witten's formulas
were confirmed a few years later by Jeffrey-Kirwan \cite{je:in}. Their
main result expressed the intersection pairings in terms
of iterated residues; the equivalence to Witten's version (as a sum
over irreducible representations) was obtained using results of Szenes
\cite{sz:it} (see Brion-Vergne \cite{br:ar1,br:ar2} for further developments in
this direction).  In Earl-Kirwan \cite{ea:com}, these results were
used to give explicit formulas for the relations in the cohomology
ring.

The formulas for intersection pairings in \cite{wi:tw} were stated not
only for the group $G=\SU(n)$, but for arbitrary compact, simply
connected Lie groups. The aim of the present paper is to give a proof
of Witten's formula in this generality, for any collection of
conjugacy classes $\Co_j$ for which the moduli space has at most
orbifold singularities.

It was observed by Atiyah-Bott that the moduli spaces $\M(\Sig;c)$
have a natural structure as an infinite-dimensional symplectic
quotient for the gauge group action on the space of connections on
$\Sig$. Witten obtained his formulas by an application of equivariant
localization techniques to this infinite-dimensional setting.
Jeffrey-Kirwan worked with a different expression of the moduli space
$\M(\Sig;c)$ as a finite-dimensional symplectic quotient $M\qu G$,
however the symplectic manifold $M$ is both singular and non-compact. 
This led to technical difficulties, which prevented the 
generalization of this approach to other groups and holonomies. 
A second problem was that Atiyah-Bott's result for generators of the 
cohomology ring, and its subsequent generalization by 
Biswas-Raghavendra, was not established for general 
compact simple Lie groups. 

In this paper, we will calculate the intersection pairings using 
localization on smooth, compact, finite-dimensional manifolds. 
This calculation is based on a more general localization formula for
Hamiltonian $G$-spaces with {\em group-valued} moment maps $\Phi:\,M\to
G$, as introduced in \cite{al:mom}. In the case at hand,
$\M(\Sig;\Co_1,\ldots,\Co_r)$ is expressed as a symplectic quotient
$M\qu G=\Phinv(e)/G$ where $M=G^{2s}\times \Co_1\times\cdots
\times\Co_r$, with $G$ acting by conjugation and moment map
$$ \Phi(a_1,b_1,\ldots,a_s,b_s,u_1,\ldots,u_r)=\prod_{i=1}^s
[a_i,b_i]\prod_{j=1}^r u_j$$
where $[a,b]=aba^{-1}b^{-1}$ is the Lie group commutator. 
As for usual $\g^*$-moment maps, there is a localization formula
\cite{al:gr}, expressing intersection pairings of classes 
in the image of the ring homomorphism $H_G(M)\to H(M\qu G)$ 
in terms of fixed point data on $M$. Unfortunately, 
in contrast to the $\g^*$-valued case this homomorphism need not be 
surjective. For instance, it turns out that for the moduli space $\M(\Sig;c)$, 
the Atiyah-Bott classes ${\mf a}^p,\ {\mf b}^p_i$ are in the 
image, but not in general the classes ${\mf f}^p$.   
 
The missing generators can be recovered as follows. For any invariant
polynomial $p$ of degree $d>0$, transgression defines a closed
equivariant differential form $\eta^p_G\in\Om_G^{2d-1}(G)$ on $G$.
Suppose $\om^p\in\Om_G^{2d-2}(M)$ is an equivariant form with
$\d_G\om^p=\Phi^*\eta_G^p$. Then the pull-back of $\om^p$ to the level
set $\Phinv_G(e)$ is equivariantly closed, hence defines a class in
$H_G(\Phinv_G(e))\cong H(M\qu G)$.  A recent result of
Bott-Tolman-Weitsman \cite{bo:su} asserts that classes of this type,
together with the image of the Kirwan map, generate the cohomology
ring of $M\qu G$ provided {\em all} forms $\Phi^*\eta_G^p$ are
equivariantly exact.  In Section \ref{sec:witproof}
we will show that this condition holds for
the moduli space examples, by explicit construction of the forms
$\om^p$ for these cases. That is, one has canonical generators for the
cohomology ring of $\M(\Sig;\Co_1,\ldots,\Co_r)$ in full generality.

The main result of this paper is a localization formula for
intersection pairings of the Bott-Tolman-Weitsman classes on $M\qu G$
with classes in the image of the Kirwan map. In the
moduli space setting, the evaluation of the fixed point data is
fairly straightforward and immediately leads to Witten-type formulas.

The organization of the paper is as follows. In Section
\ref{sec:group} we review $G$-valued moment maps, and in Section
\ref{sec:abel} we show how to relate $G$-valued moment maps to
standard $\g^*$-moment maps (linearization) or to $T$-valued moment
maps (Abelianization).  Section \ref{sec:dh} is a review of
Duistermaat-Heckman theory for $\g^*$-valued moment maps. The somewhat
unusual perspective taken here is that Duistermaat-Heckman measures
are cohomology classes for a 'twisted' equivariant differential on
$\g^*$. We explain that the DH-distributions can be equally well defined
as distributions on $\t^*$, and give a geometric interpretation. These
ideas are put to use in Section \ref{sec:qdh}, in order to define
DH-distributions for group-valued moment maps. Again one encounters a
certain twisted differential, this time on the group $G$ rather than
on $\g^*$. We associate to any cocycle for this differential an
invariant distribution on $G$, by first Abelianizing the problem. We
find that on one hand, the resulting DH-distributions encode
intersection pairings on symplectic quotients, and on the other hand
are calculated by localization. In Section 6, we use a
similar approach to study more complicated twistings, required to deal
with the Bott-Tolman-Weitsman classes $[\om^p_{\on{red}}]$. As it
turns out, Witten's 'change of variables' \cite{wi:tw} becomes very
natural in this context. Finally, in Section \ref{sec:witproof}
we apply our localization formula to the moduli space examples.
\vskip.2in 
\newpage
\noindent{\bf Acknowledgments.} It is a pleasure to thank Anton
Alekseev, Lisa Jeffrey and Chris Woodward for many helpful discussions
and useful comments. Some of these results were obtained during a stay
at Chuo University in April 1998, and I would like to thank the
Mathematics Department, and particularly Professor Takakura, for their
hospitality.
\vskip.4in

\noindent{\bf Notation.}  Throughout this paper, $G$ will be a
compact, connected Lie group.  We fix an invariant inner product
$\cdot$ on the Lie algebra $\g$, which we will often use to identify 
$\g$ with its dual $\g^*$.  Choose a maximal torus $T$ of $G$, with Lie
algebra $\t$, and let $\mf{p}=\t^\perp$ be the unique $T$-invariant
complement of $\t$ in $\g$. The Weyl group of $(G,T)$ will be denoted 
$W=N_G(T)/T$. 

The integral lattice $\Lambda\subset\t$ is
the kernel of the exponential map $\exp:\,\t\to T$, and the (real)
weight lattice is its dual, $\Lambda^*=\{\lambda\in
\t^*|\,\l\lambda,\xi\r\in \Z\ \forall \xi\in\Lambda\}$. Weights
$\lambda\in\Lambda^*$ parametrize homomorphims $\eps_\lambda:\ T\to
\U(1)$ where
$$\eps_\lambda(\exp\xi)=e^{\tpi \lambda\cdot \xi},\ \ \xi\in\t.$$
Let $\mf{R}\subset \Lambda^*$ be the set of (real) roots.  Fix a
positive Weyl chamber $\t_+\subset\t$, and let 
$$\mf{R}_+=\{\alpha\in
\mf{R}|\,\l\alpha,\xi\r\ge 0\mbox{ for all }\xi\in\t_+\}$$ 
be the
corresponding set of positive roots.  The cardinality of the set of
positive roots will be denoted by $n_+=\#\mf{R}_+$. Recall that the
choice of $\t_+$ defines a unique $T$-invariant complex structure on
$\mf{p}$, in such a way that $\mf{R}_+$ are the weights for the
$T$-action. In particular, this defines an orientation on $\mf{p}$. 
 
Let the homogeneous space $G/T$ be equipped with the $G$-invariant
Riemannian metric and orientation induced from $\mf{p}\cong \g/\t$.
The Riemannian volume of $G/T$ is given by the formula \cite[Corollary
7.27]{be:he}
\begin{equation}\label{eq:volGrho}
\vol_{G/T}=(\prod_{\alpha\in\mf{R}_+} 2\pi\alpha\cdot\rho)^{-1},
\end{equation}
where $\rho=\hh\sum_{\alpha\in\mf{R}_+}\alpha$. More generally, if
$K\subset G$ is a connected closed subgroup containing $T$, the
homogeneous space $G/K$ carries a Riemannian metric and orientation
induced from $\k^\perp$, and the Riemannian volume $\vol_{G/K}$ is
given by a formula similar to \eqref{eq:volGrho}, but with a product
over only those $\alpha\in\mf{R}_+$ that are not roots of $K$.
 
Let the $\t_+^*\subset\t^*$ be unique Weyl chamber in $\t^*$
containing $\rho$. (Clearly, the identification of $\t$ and $\t^*$
given by the inner product identifies $\t_+\cong \t^*_+$.)
For any $\mu\in \t_+^*$, the symplectic
volume of the (co-)adjoint orbit $G.\mu$ is given by 
\begin{equation}\label{eq:volGmu}
\Vol(G.\mu)=(\prod_{\alpha\cdot\mu>0}2\pi \alpha\cdot\mu) 
\ \vol_{G/G_\mu},
\end{equation}
where $\vol_{G/G_\mu}$ is the Riemannian volume corresponding to the
metric on $\g_\mu^\perp$ induced from $\g$. 

For any $\lambda\in \Lambda^*_+:=\Lambda^*\cap\t_+$ we denote by
$V_\lambda$ the irreducible representation of highest weight
$\lambda$, and by $\chi_\lambda$ its character.
\newpage

\section{Group-valued moment maps}\label{sec:group}
\subsection{q-Hamiltonian $G$-spaces}
In this Section we recall the concept of a group-valued moment map
introduced in \cite{al:mom}.

Recall that for any $G$-action on a manifold $M$, the equivariant 
cohomology $H_G(M)$ may be computed using the Cartan complex
$(\Om_G(M),\d_G)$, where $\Om_G(M)=(\on{Pol}(\g)\otimes\Om(M))^G$ 
is the algebra of $G$-equivariant polynomial maps $\g\to \Om(M)$ 
and $\d_G$ is the equivariant differential, $(\d_G\beta)(\xi)=
(\d-\iota(\xi_M))\beta(\xi)$. See Appendix \ref{app:carmod}.  
Letting $G$ act on itself by conjugation, 
we define an equivariant 3-form $\eta_G\in\Om^3_G(G)$ by 
$$ \eta_G(\xi)=\f{1}{12}\theta^L\cdot [\theta^L,\theta^L]-
\f{1}{2}(\theta^L+\theta^R)\cdot\xi$$
where $\theta^L,\theta^R\in\Om^1(G)\otimes\g$ are the 
left-invariant, right-invariant Maurer-Cartan forms. 
It is easily verified that $\eta_G$ is an equivariant cocycle, i.e.
$\d_G\eta_G=0$.  A {\em Hamiltonian $G$-space with group-valued moment map}
(in short, a {\em q-Hamiltonian $G$-space}) is a triple $(M,\om,\Phi)$
consisting of a $G$-manifold, an equivariant map $\Phi:\,M\to G$, and
an invariant 2-form $\om\in\Om^2(M)$ satisfying the {\em moment map 
condition}
$$ \d_G\om=\Phi^*\eta_G$$
and the {\em minimal degeneracy condition}
$$\ker(\om_x)=
\{\xi_M(x)|\,\Ad_{\Phi(x)}\xi=-\xi\},\ \ x\in M.$$  
Sometimes we will omit the minimal degeneracy condition, in which case
we refer to $(M,\om,\Phi)$ as a {\em degenerate} q-Hamiltonian
$G$-space. We list some examples of non-degenerate q-Hamiltonian 
$G$-spaces, with references for further details: 
\begin{enumerate}
\item Conjugacy classes $\Co\subset G$, with moment map the 
inclusion \cite{al:mom}, 
\item $G^2$, with $G$ acting by conjugation on each factor and 
moment map 
$(a,b)\mapsto aba^{-1}b^{-1}$ the Lie group commutator \cite{al:mom}, 
\item 
For any symmetric space $X=G/K$ of $G$, the space $X^2$, with moment  
the product of the natural inclusions $X\to G$ \cite{al:du}, 
\item even dimensional spheres $S^{2n}$, viewed as compactifications of 
a ball $B\subset \C^n$, and with $G=\U(n)$-action induced from the 
defining representation on $\C^n$ (see \cite{al:du,hu:re1} for $n=2$, the 
generalization to higher rank was recently obtained by 
Hurtubise-Jeffrey-Sjamaar \cite{hu:imp}).  
\end{enumerate}

There is a symplectic reduction procedure for q-Hamiltonian
$G$-spaces, similar to the usual Marsden-Weinstein reduction for
$\g^*$-valued moment maps: If the group unit $e\in G$ is a regular
value of the moment map, then $G$ acts locally freely on the level set
$\Phinv(e)$, the pull-back of $\om$ is $G$-basic, and the induced
2-form $\om_{\on{red}}$ on $\Phinv(e)/G$ is {\em symplectic}.  We will
refer to
$$ M\qu G=\Phinv(e)/G$$
as the {\em symplectic quotient} of $(M,\om,\Phi)$. More generally, 
one defines symplectic quotients $M_g=\Phinv(g)/G_g$ 
(where $G_g$ is the centralizer of $g$) at other regular values 
of the moment map. 

It was shown in \cite{al:mom} that the moduli space
\eqref{eq:modulispace} of flat $G$-bundles 
may be written as a symplectic quotient 
$$
\M(\Sig;\Co_1,\ldots,\Co_r)=M\qu G$$ 
of a q-Hamiltonian $G$-space $(M,\om,\Phi)$. Here
\begin{equation}\label{eq:spaceM}
 M=G^{2s}\times \Co_1\times\cdots \times\Co_r,
\end{equation}
with $G$ acting by conjugation on each factor, and moment map 
\begin{equation}\label{eq:mommap}
 \Phi(a_1,b_1,\ldots,a_s,b_s,u_1,\ldots,u_r)=\prod_{i=1}^s [a_i,b_i] 
\prod_{j=1}^r u_j.
\end{equation}
Here $[a,b]\equiv aba^{-1}b^{-1}$ is the group commutator. The 2-form 
$\om$ on $M$ is given by an explicit, but somewhat complicated formula 
spelled out in \cite{al:mom} (see also Section \ref{sec:witproof} below).  

In general, q-Hamiltonian $G$-spaces need not be orientable -- 
a counterexample is $M=\R P(2)$ as a conjugacy class for 
$G=\on{SO}(3)$. However, this problem does not arise if $G$ is simply 
connected, or more generally if the half-sum of positive roots, $\rho$, is 
in the weight lattice $\Lambda^*$. For any differential form 
$\beta\in\Om(M)$, let $\beta^{[k]}$ denotes its component in $\Om^k(M)$.

\begin{lemma}\cite{al:du}\label{lem:ori}
Suppose $\rho$ is a weight of $G$. 
Let $(M,\om,\Phi)$ be a q-Hamiltonian $G$-space. Then $M$ 
carries a canonical volume form $\Gamma$ with the property
\begin{equation}\label{eq:volume}
 (\exp\om)^{[\dim M]}=\f{\Phi^*\chi_\rho}{\dim V_\rho}\ \Gamma.
\end{equation}
In particular, $M$ is orientable. 
\end{lemma}

\begin{remark}
Suppose that $M$ is connected and that $\Phi^*\chi_\rho$ does not
identically vanish on $M$. Then \eqref{eq:volume} may be used as a
definition of the volume form $\Gamma$.  In particular, this is the
case if $\Phi^{-1}(e)\not=\emptyset$.  Note that for all $x\in
\Phinv(e)$, the 2-form $\om_x$ on $T_xM$ is non-degenerate, and that
the symplectic orientation on $T_xM$ coincides with the orientation
given by $\Gamma$.
\end{remark}

\begin{example}\label{ex:conclass}
For a conjugacy class $\Co=G/G_g$, the volume form may be described
(up to sign) as the Riemannian volume form for the homogeneous space
$G/G_g$, times $|\det_{\g_g^\perp}(\Ad_g-1)|^{1/2}$.  Suppose 
$G$ is simple and simply connected, and let
$\Alc\subset\t$ be the fundamental alcove, i.e. the subset of $\t_+$
cut out by the inequality $\alpha_{\on{max}}\cdot\mu\le 1$ where
$\alpha_{\on{max}}$ is the highest root. Recall that $\Alc$
parametrizes the set of conjugacy classes in $G$, in the sense that
every conjugacy class contains a unique element $\exp(\mu)$ with
$\mu\in \Alc$.  The volume of the conjugacy class $G.\exp\mu$ is given
by the formula \cite{al:du},
\begin{equation}\label{eq:covol}
 \Vol(G.\exp\mu)=\big(\prod_{\alpha\in\mf{R}_+,\,\,
\alpha\cdot\mu\not\in\Z} 2\sin \pi\alpha\cdot\mu\big)\,\vol_{G/G_{\exp\mu}}
\end{equation}
compare with \eqref{eq:volGmu}.  The orientation on the conjugacy
class $G.\exp\mu$ differs from the orientation of the homogeneous
space $G/G_{\exp\mu}$ by a sign, $(-1)^{2\rho_K\cdot\mu}$, where
$2\rho_K$ is the sum of positive roots of $K=G_g$, i.e. the sum over all 
those roots $\alpha\in\mf{R}_+$ with $\alpha\cdot\mu\in\Z$. 
\end{example}

\subsection{Bott-Tolman-Weitsman theorem}\label{subsec:btw}
Suppose $(M,\om,\Phi)$ is a compact, connected, q-Hamiltonian
$G$-space, and that $e$ is a regular value of the moment
map. Pull-back to the identity level set defines a {\em Kirwan map}
$$ H^\bullet_G(M)\to H^\bullet_G(\Phi^{-1}(e))= 
H^\bullet(M\qu G).$$
In contrast to ordinary Hamiltonian $G$-spaces, this map is not
onto, in general. 

\begin{example}
Let $G=\SU(2)$, $\Co=\{c\}$ the conjugacy class consisting of the
non-trivial central element. For $s\ge 2$ consider $M=G^{2s}\times
\Co$, with moment map $\Phi$ as in \eqref{eq:mommap}.  It is easy to
see that away from $\Phinv(c)$, the $G$-action has constant
stabilizer equal to the center $Z(G)$. In particular, $G/Z(G)$ acts freely on
$\Phinv(e)$, and therefore $M\qu G$ is a smooth symplectic manifold of
dimension $(2s-2)\dim G>0$. The symplectic form defines a non-trivial 
class in  $H^2(M\qu G)$. On
the other hand $H^2_G(M)=0$. This shows that the Kirwan map cannot be
surjective in this example.
\end{example}

We now explain how to obtain a set of generators of $H(M\qu G)$ in the
general case. As it turns out, it is necessary to take the topology of
$G$ (as a target of the moment map) into account. Recall that
generators of $(\wedge\g^*)^G\cong H(G)$ are obtained as images of the
{\em transgression map} $\on{Pol}^\bullet(\g)^G\to
(\wedge^{2\bullet-1}\g^*)^G$. The transgression construction can be
made equivariant for the conjugation action: That is, there is a
canonical linear map
$$ \on{Pol}^\bullet(\g)^G\to \Om^{2\bullet-1}_G(G),\ p\mapsto \eta^p_G$$
such that the forms $\eta^p_G$ are closed, and their classes
$[\eta^p_G]$ generate $H_G(G)$ as an algebra over $H_G(\pt)$. (In
fact, one already obtains a set of generators if one restricts to the
subspace $P\subset \on{Pol}^\bullet(\g)^G$ spanned by primitive
generators for the algebra $\on{Pol}^\bullet(\g)^G$.) 

An explicit formula for the forms $\eta^p_G$ is worked out in
Jeffrey's paper \cite{je:gr}. 
For any $p\in \on{Pol}(\g)$ let
$\xi\mapsto p'(\xi)\in \g$ be its gradient, defined by
$ p'(\xi)\cdot \zeta=\f{d}{d t}\Big|_{t=0}\,p(\xi+t\zeta).$
Then
\begin{equation}\label{eq:jef}
 \eta^p_G(\xi)=-\theta^L\cdot \int_0^1 \d t \ p'\Big(
(1-t)\xi+\Ad_{g^{-1}}(\xi)-\f{t(1-t)}{2}[\theta^L,\theta^L]
\Big). 
\end{equation}
We will review this derivation of this formula in the Appendix, 
Section \ref{subsec:pt}.
\begin{remarks}\label{rem:eta}
\begin{enumerate}
\item\label{it:torus}
Note that if $G=T$, the formula for the
equivariant forms $\eta^p_G$ simplifies to
$$ \eta^p_T(\xi)=-\theta_T\cdot p'(\xi),$$
where $\theta_T$ is the Maurer-Cartan form for $T$. 
\item 
For the quadratic polynomial $p(\xi)=\hh \xi\cdot\xi$, the form $\eta_G^p$
coincides with the form $\eta_G$ considered above. 
\end{enumerate}
\end{remarks}
Suppose that $(M,\om,\Phi)$ is a q-Hamiltonian $G$-space. 
We will refer to equivariant forms $\om^p\in \Om^{2d-2}_G(M)$
with 
\begin{equation}\label{eq:highform} 
\d_G\om^{p}=\Phi^*\eta^{p}_G.
\end{equation}
as {\em higher q-Hamiltonian forms}.  Since $\eta^p_G$ has odd degree,
its pull-back to the group unit $e\in G$ vanishes and therefore
$\iota_{\Phinv(e)}^*\om^p$ is closed. Hence
$[\iota_{\Phinv(e)}^*\om^p]\in H_G(\Phinv(e))$ is defined, and
descends to an ordinary cohomology class $[\om^p_{\on{red}}]$ on the
symplectic quotient.  After choosing a principal connection on
$\Phinv(e)$, Cartan's theorem (see Appendix \ref{app:car}) yields
$\om^p_{\on{red}}\in \Om^{2d-2}(M\qu G)$ as a differential form. The
following is a reformulation of a result of \cite{bo:su}.

\begin{theorem}[Bott-Tolman-Weitsman] 
Suppose $G$ is simply connected.  Let
$(M,\om,\Phi)$ be a compact connected q-Hamiltonian $G$-space, with $e$ a
regular value of the moment map. Assume that for {\em all} $p\in
\on{Pol}(\g)^G$ there exists a form $\om^p\in \Om_G(M)$ satisfying 
\eqref{eq:highform}. Then the
classes $[\om^p_{\on{red}}]$, together with the image of the Kirwan
map, generate the cohomology ring of $M\qu G$. 
\end{theorem}
Note that $\Phi^*\eta^p_G$ is exact for all $p\in \on{Pol}(\g)^G$, if and
only if it is exact for some set of generators of the ring
$\on{Pol}(\g)^G$. That is, it is enough to consider forms $\om^p$ for such
a set of generators. 

\begin{corollary}
If $G=\SU(2)$, the image of the Kirwan
map together with the reduced symplectic form generate the
cohomology ring of $M\qu G$. 
\end{corollary}

\begin{remark}
The condition \eqref{eq:highform} means that the pair $(\om^p,\eta^p_G)$ 
defines a cocycle for the relative equivariant de Rham complex 
$\Om^\bullet_G(\Phi)=\Om^{\bullet-1}_G(M)\oplus \Om^\bullet_G(G)$. 
The class $[\om^p_{\on{red}}]$ depends only on the relative cohomology 
class $[(\om^p,\eta^p_G)]$.
\end{remark}

In Section \ref{sec:witproof}, we will give
an explicit construction of the higher q-Hamiltonian forms $\om^p$ for
the moduli space example \eqref{eq:spaceM}.

\section{Linearization and Abelianization}
\label{sec:abel}
In this Section, we will describe two methods of relating q-Hamiltonian 
$G$-spaces to more standard Hamiltonian spaces: Linearization (replacing 
the target $G$ of the moment map with $\g\cong\g^*$), and Abelianization
(replacing $G$ by the maximal torus $T$). 

\subsection{The linearization of a q-Hamiltonian $G$-space}
\label{subsec:log}\label{subsec:lin}
For any invariant polynomial $p\in \on{Pol}^\bullet(\g)^G$ let $\eta^p_\g\in \Om^{2\bullet-1}_G(\g)$ 
denote the (exact) equivariant form, 
$$ \eta^p_\g(\xi)=-\d\l\cdot,p'(\xi)\r.$$
Thus, $\eta^p_\g$ is a linearized version of the form $\eta^p_G$. Let
\begin{equation}\label{eq:varpip}
\varpi^p=\on{h}(\exp^*\eta^p_G-\eta^p_\g)\in\Om_G^{2\bullet-2}(\g),
\end{equation}
where $\on{h}:\,\Om_G^\bullet(\g)\to \Om_G^{\bullet-1}(\g)$ is the $G$-equivariant homotopy operator 
(Section \ref{app:hom})
for the vector space $\g$. 
For $p(\xi)=\hh ||\xi||^2$ we omit the superscript $p$, 
writing $\eta_\g=-\l\cdot,\xi\r$ and 
\begin{equation}\label{eq:varpi}
\varpi=\on{h}(\exp^*\eta_G-\eta_\g)\in\Om_G^{2}(\g).
\end{equation}
Since  $(\exp^*\eta_G-\eta_\g)^{[1]}=0$, 
one has $\varpi^{[0]}=0$, i.e. the equivariant 2-form 
$\varpi$ is an ordinary invariant 2-form.

Suppose now that $(M,\om,\Phi)$ is a q-Hamiltonian $G$-space, possibly
degenerate. Let $V\subset G$ be an invariant open neighborhood of
$e\in G$, given as the diffeomorphic image of an invariant open
neighborhood $V_0$ of $0\in\g$ under the exponential map.  Let
$\on{log}:\,V_0\to V$ be the inverse map. Replacing $M$ with
$\Phinv(V)$ if necessary, assume that $\Phi$ takes values in $V$, and
let
$$\Phi_0=\on{log}(\Phi),\ \ \om_0=\om-\Phi^*\varpi.$$
Clearly, $\d_G\om=\Phi^*\eta_G$ implies $\d_G\om_0=\Phi_0^*\eta_\g$,
which is the usual moment map condition for a Hamiltonian $G$-space.
We will refer to the Hamiltonian $G$-space $(M,\om_0,\Phi_0)$ as the
{\em linearization} of the q-Hamiltonian space $(M,\om,\Phi)$.
Given a higher q-Hamiltonian form $\om^p$ on $M$, the form 
$$ \om_0^p=\om^p-\Phi_0^*\varpi^p$$
has the property $\d_G\om_0^p=\Phi_0^*\eta^p_\g$. 

\begin{remarks}
\begin{enumerate}
\item
If $\om$ satisfies the minimal degeneracy
condition, then the
2-form $\om_0$ is symplectic \cite[Proposition 3.4]{al:mom}. 
Let $\Gamma_0=(\exp\om_0)^{[\dim M]}$
denote the symplectic (Liouville) volume form on $M_0$, and $\Gamma$
the volume form on $M$. Then $\Gamma=\Phi_0^*J^{1/2}\,\Gamma_0$, where
$J^{1/2}\in C^\infty(\g)$ is the unique smooth square root with
$J^{1/2}(0)=1$. See \cite[Section 3.6]{al:du}.
\item
As a typical application of the linearization construction, the
symplectic reduction theorem for q-Hamiltonian $G$-spaces 
follows directly from the usual
Hamiltonian setting, together with the fact that $\varpi$ vanishes at
$\mu=0$.
\end{enumerate}
\end{remarks}

\subsection{The Abelianization of a q-Hamiltonian $G$-space}
\label{subsec:abel}
The theory of $G$-valued moment maps is substantially different from
the theory of ordinary $\g^*$-valued moment maps only if the group $G$
is non-Abelian. In this Section we introduce an Abelianization
procedure for q-Hamiltonian $G$-spaces, replacing the group $G$ by its
maximal torus $T$. For any $G$-manifold 
$M$, we denote by $\kappa_T:\,\Om_G(M)\to \Om_T(M)$ the natural map 
restricting the action. 

Recall that $\mf{p}=\t^\perp$, and consider the two maps
$$ 
\xymatrix{&T\times \mf{p}\ar[ld]_{\pi_T} \ar[rd]^F&\\
T&&G}
$$
where $\pi_T(t,\mu)=t$ and $F(t,\mu)=t\exp(\mu)$. 
For any $p\in \on{Pol}^\bullet(\g)^G$, define 
a $T$-equivariant form 
\begin{equation}\label{eq:gammap}
 \gamma^p=\on{h}\Big(F^*\kappa_T(\eta_G^p)\Big),\ \ 
\in \Om_T^{2\bullet-2}(T\times\mf{p})
\end{equation}
where $\on{h}$ is the $T$-equivariant homotopy operator 
(Section \ref{app:hom}) for
$T\times\mf{p}\to T$. 
\begin{proposition}
The pull-back of the forms $\gamma^p$ to $T\times\{0\}\subset T\times\mf{p}$
vanishes. One has, 
$$ \d_T\gamma^p=F^*\kappa_T(\eta_G^p)-\pi_T^*\eta_T^p,\ \ 
\xi\in\t.$$
\end{proposition}

\begin{proof}
Let $\iota_T:\,T\to T\times \mf{p}$ denote the inclusion.  
The equation for $\d_T\gamma^p$ follows from the property
$\d_T\on{h}+\on{h}\d_T=\on{id}-\pi_T^*\iota_T^*$ of the homotopy
operator, since $\iota_T^*(F^*\kappa_T(\eta_G^p))=\eta_T^p$.  Note
that one can also write
$\gamma^p=\on{h}\Big(F^*\kappa_T(\eta_G^p)-\pi_T^*\eta_T^p\Big)$ since
$\on{h}\circ \pi_T^*=0$. The pull-back of
$F^*\kappa_T(\eta_G^p)-\eta_T^p$ to $T$ is zero, hence the same is
true for $\gamma^p$.
\end{proof}

Consider in particular the equivariant 2-form 
$\gamma\in \Om^2_T(T\times\mf{p})$ 
corresponding to the quadratic polynomial 
$p(\xi)=\hh ||\xi||^2$. It turns out that $\gamma^{[0]}=0$ 
so that $\gamma$ is in ordinary 2-form: 

\begin{proposition}
The equivariant 2-form $\gamma$ is a $T$-invariant 2-form, given by 
the formula, 
$$ \gamma=\varpi_{\mf{p}}+\hh \theta_T\cdot 
\exp_{\mf{p}}^*\theta^R\in \Om^2(T\times\mf{p}),$$
where $\varpi_{\mf{p}}$ is the pull-back of the form 
$\varpi\in\Om^2(\g)^G$ to $\mf{p}$, and $\exp_{\mf{p}}$ is the restriction 
of the exponential map. The pull-back of $\gamma$ to $T$ vanishes.  
\end{proposition}
\begin{proof}
We use the following formula for the pull-back of the form 
$\eta_G$ under group multiplication: 
$$ \on{Mult}_G^*\eta_G=\pi_T^*\eta_G +\pr_2^*\eta_G+\hh d_G 
(\pi_T^*\theta^L\cdot \pr_2^*\theta^R).$$
Here $\pr_i$ denote the projections from $G\times G$ to the two 
factors. The map $F$ can be written as a composition of the map 
$T\times\mf{p}\to G\times G,\,(t,\mu)\mapsto (t,\exp_{\mf{p}}(\mu)$ 
followed by group multiplication. Hence, 
$$ F^*\kappa_T(\eta_G)=\eta_T+\exp_{\mf{p}}^*\kappa_T(\eta_G)+\hh \d_T(\theta_T \cdot
\exp_{\mf{p}}^*\theta^R).$$
Now apply the homotopy operator $\on{h}$ for $T\times\mf{p}$. 
Clearly, $\on{h}\eta_T=0$. Furthermore, 
$$ \on{h}\exp_{\mf{p}}^*\kappa_T(\eta_G)=\varpi_{\mf{p}}$$ 
since pull-back from $\g$ to $\mf{p}$ intertwines the 
homotopy operators for the two vector spaces. Finally, 
$$ \on{h}\d_T(\theta_T \cdot
\exp_{\mf{p}}^*\theta^R)=\theta_T \cdot
\exp_{\mf{p}}^*\theta^R$$
since $\theta_T \cdot
\exp_{\mf{p}}^*\theta^R$ pulls back to $0$ on $T\subset T\times\mf{p}$, 
and since 
$$ \on{h} (\theta_T \cdot
\exp_{\mf{p}}^*\theta^R)=-\theta^T\cdot\mu=0.$$
\end{proof}

We may use the form $\gamma$ to turn any q-Hamiltonian $G$-space
$(M,\om,\Phi)$ into a (degenerate) q-Hamiltonian $T$-space, at least
on a neighborhood of $\Phinv(T)$. 

For $\eps>0$ let $B_\eps(0)\subset \mf{p}$ denote the open ball of 
radius $\eps$. Choose $\eps$ sufficiently small, so that the  
map $F$ restricts to a 
a diffeomorphism $T\times B_\eps(0)$ 
onto an open subset $U\subset G$. Let $\pi:\,U\to T$ be the 
projection corresponding to $\pi_T:\,T\times\mf{p}\to T$, 
and use $F$ to view $\gamma^p$ as forms on $U$. Thus, 
$$ \d_T\gamma^p=\kappa_T(\eta_G^p)|_U-\pi^*\eta_T^p.$$
The pre-image 
$$ N:=\Phinv(U)$$
is a $T$-invariant open neighborhood of $\Phinv(T)$ in $M$. 
Define an invariant 2-form $\om_N\in \Om^2(N)$ and an equivariant map  
$\Phi_N\in C^\infty(N,\, T)$ by 
$$ \om_N=\om|_N-\Phi^*\gamma,\ \ 
\Phi_N=\pi\circ \Phi|_N. $$
More generally, if $\om^p$ is a higher q-Hamiltonian form for 
an invariant polynomial $p$, i.e. $\om^p$ 
is a primitive for $\Phi^*\eta^p_G$, we define 
$$ \om_N^p=\kappa_T(\om^p)|_N-\Phi^*\gamma^p.$$
\begin{proposition}\label{prop:abel}
The triple $(N,\om_N,\Phi_N)$ is a (degenerate) q-Hamiltonian $T$-space. 
If $p\in \on{Pol}(\g)^G$ is an invariant polynomial 
and $\om^p$ a corresponding higher q-Hamiltonian form for 
$(M,\om,\Phi)$, then $\om^p_N$ is a higher q-Hamiltonian form 
for $(N,\om_N,\Phi_N)$. 
\end{proposition}

\begin{proof}
The equation $\d_T \om_N^p=\Phi^*\eta^p_T$
follows directly from $\d_G\om^p=\Phi^*\eta_G^p$ and the 
property of the form $\gamma^p$. 
\end{proof}
Notice that the 2-form $\om_N$ is closed, but is usually not symplectic 
even if $\om$ is minimally degenerate.

\section{DH-distributions for Hamiltonian $G$-spaces}\label{sec:dh}
For any manifold $M$, let $\ca{C}^\bullet(M)$ denote the complex of de
Rham {\em currents} on $M$. Thus $\ca{C}^k(M)$ is the topological dual
space of $\Om^{n-k}(M)_{\on{comp}}$, and the differential is defined
by duality. If $M$ is oriented, the natural inclusion
$\Om^\bullet(M)\to \ca{C}^\bullet(M)$ is a quasi-isomorphism.  For any
vector field $X$ on $M$, the operators of Lie derivative and
contraction $L_X,\,\iota_X$ extend to the complex $\ca{C}^\bullet(M)$.
Given a $G$-action on $M$, the space $\ca{C}^\bullet(M)$ is a
$G$-differential space in the sense of \cite{ku:eq}, and therefore a 
differential space
$$\ca{C}_G^\bullet(M)=(\on{Pol}(\g)\otimes \ca{C}^\bullet(M))^G$$
of equivariant currents is defined. 
 
In the present Section, we will avoid identifications of the Lie
algebra with its dual.  Suppose $(M,\om,\Phi)$ is a Hamiltonian
$G$-space, i.e. with an ordinary moment map $\Phi:\,M\to \g^*$. 
For any equivariant cocycle $\beta\in \Om_G(M)$, the expression $\beta
e^{\tpi \om}\in\Om_G(M)$ is a cocycle for the twisted differential,
$\d_G-\tpi \Phi^*\eta_\g$.  (Note that the equivariant 3-form
$\eta_\g(\xi)=-\d\l\cdot,\xi\r$ is invariantly defined as a form on
$\g^*$.)  Hence, its push-forward under the moment map is an
equivariant current on $\g^*$, closed under the differential
$\d_G'=\d_G-\tpi \eta_\g$. We will show in this Section that the
cohomology space of this differential is naturally identified with
invariant distributions. The distributions associated with the 
current $\Phi_*(\beta e^{\tpi \om})$ will be called Duistermaat-Heckman 
distributions. 

\subsection{Equivariant currents on $\g^*$}
Even though $\eta_\g$ is equivariantly 
exact, it still defines a non-trivial twisting of the cohomology, as the 
following result shows. 

\begin{proposition}\label{prop:twisted}
Let $\d_G'=\d_G-\tpi \eta_\g$.
For any $G$-invariant open subset $V\subset \g^*$, the map
\begin{equation}\label{eq:qis}
\ca{C}_G(V)\to \D'(V)^G,\ \phi\mapsto \n=\sum_I
\big(\f{1}{\tpi}\f{\p}{\p \mu}\big)^I \phi^{[\dim
\g]}_I.\end{equation}
vanishes on $\d_G'$-coboundaries, and descends to an isomorphism, 
$$ H(\ca{C}_G(V),\d_G')\to \D'(V)^G.$$
If $\phi\in \ca{C}_G(V)$ is compactly supported, the 
associated distribution $\n$ is given in terms of its Fourier 
transform by 
\begin{equation}\label{eq:integral1}
 \l\n,e^{-\tpi\l\cdot,\xi\r}\r=
\l \phi(\xi),\,e^{-\tpi\l\cdot,\xi\r}\r.
\end{equation}
The map \eqref{eq:qis} restricts to a map for smooth equivariant forms, 
$\Om_G(V)\to \Om^{\dim\g}(V)^G$, which descends to 
to an isomorphism $H(\Om_G(V),\d_G')\to (\Om^{\dim\g}(V))^G$.
\end{proposition}
\begin{proof}
The formula \eqref{eq:integral1} for compactly supported $\phi$ follows from 
$$ \l \n,\,e^{-\tpi\l\cdot,\xi\r}\r=
\big\l\sum_I
\big(\f{1}{\tpi}\f{\p}{\p \mu}\big)^I \phi_I,\,e^{-\tpi\l\cdot,\xi\r}
\big\r
=\sum_I \xi^I \l \phi_I, e^{-\tpi\l\cdot,\xi\r}\r.$$
Note that the right hand side of \eqref{eq:integral1} may also be viewed 
as the integral (push-forward to a point) of the current 
$\phi(\xi) e^{-\tpi\l\cdot,\xi\r}$. 
We next show that the distribution $\n$ corresponding to a coboundary
$\d_G'\psi$ is zero. Using a partition of unity, it suffices to prove
this if $\psi$ has compact support. But in this case, the integral of
$(\d_G'\psi)(\xi) e^{-\tpi\l\cdot,\xi\r}
=\d_\xi(\psi(\xi)e^{-\tpi\l\cdot,\xi\r})
$ vanishes by Stokes' theorem.  It remains to show that
the induced map in $\d_G'$-cohomology is an isomorphism. To this
end, we view the space $\ca{C}(V)$ of currents as a $G$-differential space, 
with the standard $G$-action and the standard differential, but 
with the contraction operators
$$\iota_\xi'=\iota(\xi_{\g^*})+\tpi \eta_\g(\xi).$$
>From this perspective, the
$\d_G'$-cohomology of $\ca{C}_G(V)$ is just the equivariant cohomology
of the $G$-differential space $\ca{C}_G(V)$.  
Let $P \in C^\infty(\g^*,\wedge^2\g^*)$ be the Kirillov-Poisson bivector 
field on $\g^*$. If $f_{ab}^c$ are the structure constants 
in a given basis $e_a$ of $\g$, with dual basis $e^a$, 
and $\mu_a$ the associated coordinates, 
$$ P=\hh \sum_{abc} f_{ab}^c \mu_c \f{\p}{\p\mu_a}\f{\p}{\p\mu_b}.$$
Let $\iota(P):\,\ca{C}^\bullet(V)\to \ca{C}^{\bullet-2}(V)$ 
denote the operator of contraction by $P$. Then conjugation 
by $\exp(\f{1}{\tpi}\iota(P))$ is a $G$-equivariant automorphism of 
$\ca{C}(V)$ which simplifies the contraction operators: Indeed, 
$$ [\iota(P),\,\iota(\xi_{\g^*})]=0,\ \ 
[\iota(P),\,\eta_\g(\xi)]=-\iota(\xi_{\g^*}),$$
and therefore 
$$ \iota_\xi'':=
\Ad\Big(\exp(\f{1}{\tpi}\iota(P))\Big)\iota_\xi'=\tpi\eta_\g(\xi)
=\tpi \d\l\cdot ,\xi\r
.$$
Clearly, the horizontal subspace of $\ca{C}(V)$ with respect to
$\iota_\xi''$ is the space of currents of top degree, that is 
$\D'(V)$, and so the basic subspace is $\D'(V)^G$. It is easy to see that the conjugated differential 
$\d''=\Ad\big(\exp(\f{1}{\tpi}\iota(P))\big)\d$ vanishes on $\D'(V)^G$: 
One one hand, it preserves the basic subspace, on the other hand, it 
changes parity while everything in the basic subspace has fixed 
parity given by $\dim\g$. This shows that there is an isomorphism
$$
\D'(V)^G\to \ca{C}(V)_{\on{basic}},\ \ \lambda\mapsto 
\exp(\f{1}{\tpi}\iota(P))\lambda
$$
taking values in cocycles. Viewing $\ca{C}(V)_{\on{basic}}$ as a
subspace of $\ca{C}_G(V)$, it is obvious that its composition with the
map $\ca{C}_G(V)\to \D'(V)$ is the identity. To complete the proof, it
suffices to show that the inclusion of the basic subcomplex
$\ca{C}(V)_{\on{basic}}\to \ca{C}_G(V)$ induces an isomorphism in
cohomology. By a generalized version of Cartan's theorem, due to
Guillemin-Sternberg \cite{gu:su}, this will be the case if there
exists a $G$-equivariant linear map 
$\theta:\,\g^*\to \on{End}^{\on{odd}}(\ca{C}(V))$ such that $[\iota_\xi',\theta(\mu)]=\l\mu,\xi\r\,\on{Id}$ 
and such that the operators $\theta(\mu)$ and $[\d,\theta(\mu)]$
generate a super-commutative subalgebra of $\on{End}(\ca{C}(V))$. 
\footnote{These conditions are a translation of the notion of a $W^*$-module, 
as introduced in \cite{gu:su}. As shown in \cite{al:cw}, one may 
in fact drop the super-commutativity condition.}. 
Indeed, such a map is given by contraction with the constant
vector fields defined by elements of $\g^*$:
$$\theta(e^a)=\f{1}{\tpi}\iota(\f{\p}{\p\mu_a}).$$
Clearly, all of the above goes through for $\ca{C}_G(V)$ replaced with 
$\Om_G(V)$. 
\end{proof}
\begin{remark}
The map \eqref{eq:qis} restricts to a map for smooth equivariant forms, 
$\Om_G(V)\to \Om^{\dim\g}(V)^G$, and the proof of Proposition 
\ref{prop:twisted} shows that, again, this map descends to 
to an isomorphism $H(\Om_G(V),\d_G')\cong \Om^{\dim\g}(V)^G$.
\end{remark}
Later we will need the following consequence of the proof of Proposition 
\ref{prop:twisted}: 
\begin{corollary}\label{cor:cutoff}
Suppose $\phi\in \ca{C}_G(V)$ is a cocycle for $\d_G'=\d_G-\tpi
\eta_\g$, and $\O\subset V$ is a coadjoint orbit. Then there 
exists a $\d_G'$-cocycle $\ti{\phi}$ with compact support in any 
given neighborhood of $\O$, such that $\ti{\phi}=\phi$ on a smaller
neighborhood of $\O$. If $\phi$
is smooth near $\O$, one can take $\ti{\phi}$ to be smooth. 
\end{corollary}
\begin{proof}
Let $\n\in\D'(V)^G$ be the distribution defined by $\phi$. Let 
$V_1,V_2$ be open neighborhoods of $\O$, with $\ol{V_1}\subset V$ 
and $\ol{V_2}\subset V_1$. Let $\chi\in C^\infty(V)^G_{\on{comp}}$
with $\chi|_{V_1}=1$, and $\chi_1 \in C^\infty(V_1)^G_{\on{comp}}$
with $\chi_1|_{V_2}=1$. Let
$\n_1=\chi\,\n$, and let $\phi_1=\exp(\f{1}{\tpi}\iota(P))\n_1$ as in
the proof of Proposition \ref{prop:twisted}. Since $\n_1=\n$ on
$V_1$, it follows that the restrictions of $\phi,\phi_1$ to $V_1$ are
$\d_G'$-cohomologous: $\phi|_{V_1}=\phi_1|_{V_1}+\d_G'(\psi)$.
Define $\ti{\phi}$ by $\phi_1|_{V_1}+\d_G'(\chi_1
\psi)$ on $V_1$, and equal to $\phi_1$ outside $V_1$.
\end{proof}

\subsection{DH-distributions}\label{subsec:DH1}
Suppose $(M,\om,\Phi)$ is an oriented Hamiltonian $G$-space, with a
possibly degenerate 2-form $\om$.  We assume that the moment map
$\Phi:\,M\to \g^*$ takes values in an invariant open subset $V\subset
\g^*$, and is proper as a map into $V$.  For any equivariant cocycle
$\beta\in\Om_G(M)$, the push-forward
$$ \Phi_*(\beta e^{\tpi \om})\in \ca{C}_G(V)$$ 
is a well-defined $\d_G'$-cocycle. We define the Duistermaat-Heckman 
distribution $\n^\beta\in \D'(V)^G$ to be its image under the map 
\eqref{eq:qis}, that is, 
$$ \n^\beta=\sum_I (\f{1}{\tpi}\f{\p}{\p\mu})^I\,\Phi_*(\beta
e^{\tpi \om})^{[\dim M]}_I.$$
For $\beta=1$, this simplifies to $\n=\Phi_*({e^{\tpi
\om}})^{[\dim M]}$, which (up to $\tpi$ factors) is the original
definition of the Duistermaat-Heckman measure in \cite{du:on}, as a
push-forward of the Liouville measure.  The more general
DH-distributions $\n^\beta$ were introduced by Jeffrey-Kirwan
\cite{je:lo1}, in terms of their Fourier coefficients (see Proposition
\ref{prop:nb0}\eqref{it:integral} below).

We list some of the basic properties of DH-distributions. 
Given $\beta$, let us call $\mu\in \g^*$ a {\em $\beta$-regular value} of 
$\Phi$ if the differential $\d_x\Phi$ has maximal rank for all 
$x\in \Phi^{-1}(\mu)\cap \on{supp}(\beta)$. This includes of course 
the set of regular values of $\Phi$.
\begin{proposition}[Properties of Hamiltonian DH-distributions]\label{prop:nb0}
Let $(M,\om,\Phi)$ be an oriented Hamiltonian $G$-space (with possibly
degenerate 2-form $\om$), with proper moment map $\Phi:\,M\to V\subset
\g^*$. 
\begin{enumerate}
\item\label{it:locality}
The definition of $\n^\beta$ is local, in the sense that the 
restriction of $\n^\beta$ to an invariant open subset $V'\subset V$ is 
the DH-distribution associated to $\beta|_{\Phinv(V')}$. 
Its support and singular support satisfy 
$$\on{supp}(\n^\beta)\subset \Phi(\on{supp}(\beta)),\ \ 
\on{singsupp}(\n^\beta)\subset \Phi(\on{supp}(\beta)\cap (M\backslash M_*)),$$
where $M_*\subset M$ is the set of $x\in M$ such that $\d_x\Phi$ has
maximal rank.
\item \label{it:cohomology}
The map $\Om_G(M)\to \D'(V)^G,\ \beta\mapsto \n^\beta$ vanishes on 
$\d_G$-coboundaries. In particular, it descends to cohomology. 
\item \label{it:polynomial} Suppose $p\in \on{Pol}(\g)^G\subset \Om_G(M)$ is 
an invariant polynomial, viewed as an equivariant cocycle on $M$. 
Then 
$$ \n^{p\beta}=p\big(\f{1}{\tpi}\f{\p}{\p\mu}\big)\,\n^\beta.$$
\item\label{it:integral} If $\beta$ has compact support, then so does
$\n^\beta$ and its Fourier transform is given by
$$ 
\l \mf{n}^\beta,e^{-\tpi\l\cdot,\xi\r}\r=
\int_M \beta(\xi) e^{\tpi (\om-\l\Phi,\xi\r)},\ \ \xi\in\g.
$$
\item\label{it:interpret}
The distribution $\mf{n}^\beta$ is smooth near $\beta$-regular values of 
$\Phi$. If $\d_G\beta=0$, and $\mu\in\t_+^*\subset \g^*$ is a 
$\beta$-regular value, 
%
$$
\mf{n}^\beta(\mu)=(\tpi)^{\dim G_\mu/2}
\f{\vol_G}{\Vol(G.\mu)}\,\sum_j \f{1}{k_j}
\int_{\Phi^{-1}(\mu)_j/G_\mu}\beta_{\on{red}}e^{\tpi \om_{\on{red}}} 
$$
%
Here the sum is over connected components $\Phi^{-1}(\mu)_j$ of the 
level set and $k_j$ is the number of elements in a principal stabilizer for 
the $G_\mu$-action on $j$th component. Furthermore, 
$\beta_{\on{red}}$ is the image of $\beta$ under the chain map
$\Om_G(M)\to \Om(\Phi^{-1}(\mu)/G_\mu)$, 
given as pullback to the level set followed by the Cartan map for 
a principal $G_\mu$-connection on $\Phi^{-1}(\mu)$, and similarly 
for $\om_{red}$.
\end{enumerate}  
\end{proposition}
\begin{proof}
Properties \eqref{it:locality}, \eqref{it:cohomology} and 
\eqref{it:polynomial} are clear from the definition or from the 
corresponding properties of $\d_G'$-closed currents on $V$. 
The integral formula \eqref{it:integral} follows from 
$$ \l\n^\beta,e^{-\tpi\l\cdot,\xi\r}\r=
\l \Phi_*(\beta(\xi) e^{\tpi \om}),e^{-\tpi\l\cdot,\xi\r}\r=
\int_M \beta(\xi) e^{\tpi (\om-\l\Phi,\xi\r)}.$$
Proofs of part \eqref{it:interpret} may be found in \cite{je:lo1,ve:no,du:eq}.
\end{proof}

\begin{remark}
A well-known result of Kirwan \cite{ki:coh} says that if $M$ is
compact and connected, and the 2-form $\om$ is symplectic, then all
level sets $\Phi^{-1}(\mu)$ are connected. Hence, in this case
\eqref{it:interpret} directly relates pairings on
$\Phinv(\mu)/G_\mu$. More generally, this is the case if each
component of the level set $\Phinv(\mu)$ intersects the principal
stratum for the $G$-action on $M$.
\end{remark}

Since the integrand in \eqref{it:integral} is closed under
$\d_\xi=\d-\iota_{\xi_M}$, the Fourier transform of $\n^\beta$ 
for a compactly supported cocycle $\beta$ 
may be calculated by localization (Theorem \ref{app:bv}), 
as a sum over fixed point manifolds for
the vector field $\xi_M$: 
$$ 
\l \mf{n}^\beta,e^{-\tpi\l\cdot,\xi\r}\r=\sum_{F\in\ca{F}(\xi)}\int_F 
\f{\beta(\xi) e^{\tpi (\om-\l\Phi,\xi\r)}}{\Eul(\nu_F,\xi)},\ \ \xi\in\g.
$$
Combining this with the interpretation \eqref{it:interpret} of
$\n^\beta$ leads to formulas for intersection pairings on reduced
spaces in terms of fixed point data, such as the Jeffrey-Kirwan
theorem \cite{je:lo1}.

\subsection{The isomorphism $R_\g$.}
Choose orientations on $\g$ and $\t$, compatible with the orientation of 
$\mf{p}$ given by the set $\mf{R}_+$ of positive roots.  
Let $\d\vol_{\g^*}$ denote the volume form on $\g^*$, corresponding 
to the given invariant inner product on $\g$, and 
similarly $\d\vol_{\t^*}$ the volume form for the restriction of the 
inner product. 

Let $\D'(\t^*)^{W-\on{alt}}$ denote the subspace of distributions on 
$\t^*$ which {\em alternate} under the action of the Weyl group $W$, that is, 
$w_*\m=(-1)^{l(w)}\m$ where $l(w)$ is the length of $w\in W$. 
More generally, for an invariant open subset $V\subset\g^*$ we may consider 
the space $\D'(V\cap\t^*)^{W-\on{alt}}$. 

\begin{lemma}\label{lem:restriction}
There is a unique linear isomorphism 
$$R_\g:\,\D'(V)^G\to \D'(V\cap \t^*)^{W-\on{alt}},$$ 
with the following two properties: 
\begin{enumerate}
\item \label{it:lin} {\bf ($C^\infty$-linearity)} 
$R_\g(f\mf{n})=f|_{\t^*}R_\g(\n)$ for all $f\in C^\infty(V)^G$
and $\n\in \D'(V)^G$. 
\item\label{it:ft} {\bf (Fourier coefficients)}
For compactly supported $\n$, the Fourier transform of $R_\g\n$ is given 
by 
$$
\l R_\g\n,\,e^{-\tpi\l\cdot,\xi\r}\r
=(-1)^{n_+} \prod_{\alpha\in\mf{R}_+}\tpi \l\alpha,\xi\r\  
\l\n,\, e^{-\tpi\l\cdot,\xi\r}\r,\ \ \xi\in\t.
$$ 
\end{enumerate}
The image of the volume form is given by
$$R_\g(\d\vol_{\g^*})=(\prod_{\alpha\in\mf{R}_+}2\pi \alpha)\, \d\vol_{\t^*};$$
here $\alpha\in\mf{R}_+\subset \t^*$ is viewed a linear form on $\t^*$, using
the inner product.
\end{lemma}
\begin{proof}
The map on test functions 
\begin{equation}\label{eq:test}
 C^\infty(V)^G_{\on{comp}}\to
  C^\infty(V\cap \t^*)^{W-\on{alt}}_{\on{comp}},\ f\mapsto \f{\vol_{G/T}}{\#
  W} (\prod_{\alpha\in\mf{R}_+}2\pi \alpha) f|_{V\cap \t^*}
\end{equation}
is a continuous linear isomorphism, 
independent of the choice of inner product. 
Let $R_\g:\,\D'(V)^G\to \D'(V\cap \t^*)^{W-\on{alt}}$ be defined as the 
inverse of the dual map. Clearly, $R_\g$  satisfies \eqref{it:lin}. 
To establish \eqref{it:ft}, it suffices to consider
the case that $\n$ is the Liouville measure of a coadjoint
orbit $G.\mu$ for a regular element $\mu\in\t^*_+$, since linear
combinations of such measures are dense in $\D'(\g^*)^G$. That is, 
we assume that $\n$ is the unique 
invariant measure supported on $G.\mu$, with total integral equal to 
the symplectic volume 
$\Vol(G.\mu)=(\prod_{\alpha\in\mf{R}_+}2\pi \alpha\cdot\mu)
\vol_{G/T}$. In this case, 
\begin{equation}\label{eq:coadjoint}
 R_\g(\n)=\sum_{w\in W} (-1)^{l(w)} \delta_{w\mu},
\end{equation}
as one verifies on test functions. Now \eqref{it:ft} follows 
by the formula for Fourier transforms of coadjoint orbits
\cite[Chapter 7]{be:he}. Finally, the formula for 
$R_\g(\d\vol_{\g^*})$ follows from the Weyl integration formula, 
$$ \int_{\g^*} f \d\vol_{\g^*}=\f{\vol_{G/T}}{\# W}\int_{\t^*}
f|_{\t^*}\,
(\prod_{\alpha\in\mf{R}_+}2\pi \alpha)^2 \d\vol_{\t^*}.$$
\end{proof}

\begin{example}[Duistermaat-Heckman measures of coadjoint orbits]
Let $\delta_0\in \D'(\g^*)^G$ be the delta-measure supported at 
$0\in\g^*$. Then $R_\g(\delta_0)$ is expressed in terms of derivatives 
of the delta-measure on $\t^*$: 
$$ R_\g(\delta_0)=(-1)^{n_+} 
\prod_{\alpha\in\mf{R}_+}\l \alpha,\f{\p}{\p \mu}\r\ 
\delta_0.
$$
This follows from \eqref{it:ft}, by taking Fourier transforms.  
Using the characterization of $R_\g$ as the inverse dual map 
to \eqref{eq:test},  applied to $f=1$, one obtains the  identity,
\begin{equation}\label{eq:identity}
 \prod_{\alpha\in\mf{R}_+}\l \alpha,\f{\p}{\p \mu}\r\ 
\prod_{\alpha\in\mf{R}_+}2\pi \l \alpha,\cdot\r=
\f{\# W}{\vol_{G/T}}.
\end{equation}
\end{example}
\begin{proposition}
Assume $V\subset\g^*$ is an invariant open subset containing the origin. 
Suppose $\n=f\,\d\vol_{\g^*}$ with $f\in
C^\infty(V)^G$, maps to $R_\g(\n)=F\,\d\vol_{\t^*}$ with $F\in
C^\infty(V\cap\t^*)^{W-\on{alt}}$. Then 
$$ f(0)=\f{\vol_{G/T}}{\# W}
\prod_{\alpha\in\mf{R}_+}\l \alpha,\f{\p}{\p \mu}\r\Big|_{\mu=0} F.$$
\end{proposition}
\begin{proof}
This follows from \eqref{eq:identity} since 
$F= \prod_{\alpha\in\mf{R}_+}2\pi \l \alpha,\cdot\r\  f|_{\t^*}$. 
\end{proof}
By Proposition \ref{prop:twisted}, the space $\D'(V)^G$  can be thought 
of as the cohomology for the 'twisted' equivariant differential on the 
space $\ca{C}_G(V)$ of equivariant currents while 
$\D'(V\cap\t^*)$ is the cohomology of the space of equivariant 
currents on $V\cap \t^*$ (where $T$ acts trivially).
Hence, $R_\g$ may be viewed as a map of cohomology spaces, 
$H(\ca{C}_G(V),\d_G')\to H(\ca{C}_T(V\cap\t^*),\d_T')$. 
For our applications to Duistermaat-Heckman theory, it will be
important to realize this map on the 'chain level'. Let
$\pi:\,\g^*\to \t^*$ denote projection to the first factor in
$\g^*=\t^*\oplus\mf{p}^*$, and choose a representative of the
$T$-equivariant Thom form $\tau$ of the vector bundle $\pi^{-1}(V\cap
\t^*)\to V\cap \t^*$, with fiberwise compact support in the intersection 
$V\cap \pi^{-1}(V\cap \t^*)$. Recall that the pull-back of $\tau$ to $V\cap
\t^*$ represents the equivariant Euler class, hence it is
$T$-equivariantly cohomologous to $(-1)^{n_+}\prod_{\alpha\in\mf{R}_+}
\l\alpha,\cdot \r$.  Define a map
\begin{equation}\label{eq:map}
 \ca{C}_G(V)\to \ca{C}_T(V\cap\t^*)^{W-\on{alt}},\ 
\phi\mapsto \psi=(\tpi)^{n_+} \pi_*(\tau\,\kappa_T(\phi)).
\end{equation}
Clearly, \eqref{eq:map} intertwines $\d_G'$ with $\d_T'$. 
\begin{proposition}
The induced map in cohomology 
$$ H(\ca{C}_G(V),\d_G')=\D'(V)^G\to
H(\ca{C}_T(V\cap\t^*),\d_T')=\D'(V\cap\t^*)$$
takes values in
$\D'(V\cap\t^*)^{W-\on{alt}}$, and coincides with the map
$R_\g$.
\end{proposition}
\begin{proof}
Suppose the induced map in cohomology takes $\n\in \D'(V)^G$ to $\m\in
\D'(V\cap\t^*)$. We will show $\m=R_\g(\n)$.  Using a partition of
unity, we may assume that $\n$ has compact support.  As shown above,
$\n$ has a unique representative $\phi\in \ca{C}(V)_{\on{basic}}$,
given as $\phi=\exp(\f{1}{\tpi}\iota(P))\n$.  Since$\n$ has compact
support, so does $\phi$. By Fourier transform, for $\xi\in\t$,
\begin{equation}\label{eq:40}
 \l \m,e^{-\tpi\l\cdot,\xi\r}\r=
(\tpi)^{n_+}
\int_{\t^*} \pi_*(\tau(\xi)\,\phi)e^{-\tpi\l\cdot,\xi\r}=
(\tpi)^{n_+}
\int_{\g^*} \tau(\xi)\,\phi\,\, e^{-\tpi\l\cdot,\xi\r}.
\end{equation}
Since $\phi$ has compact support, we can replace the equivariant
Thom form by the cohomologous form (with non-compact support) 
$(-1)^{n_+}\prod_{\alpha\in\mf{R}_+}\l\alpha,\cdot \r\in 
\on{Pol}(\t)\subset \Om_T(\g^*)$. 
Hence \eqref{eq:40} equals 
$$ (-1)^{n_+}\prod_{\alpha\in\mf{R}_+}\tpi \l\alpha,\xi\r
\ \l\n,e^{-\tpi\l\cdot,\xi\r}\r
=\l R_\g \n,\,e^{-\tpi\l\cdot,\xi\r}\r.
$$
\end{proof}

\subsection{The distributions $\m^\beta$} \label{subsec:mbeta}
We return to the setting of Section \ref{subsec:DH1}: $(M,\om,\Phi)$
is an oriented Hamiltonian $G$-space, with proper moment map
$\Phi:\,M\to V\subset\g^*$. Let $\n^\beta\in\D'(V)$ be the invariant
distribution associated to an equivariant cocycle $\beta\in \Om_G(M)$.
We may then define a distribution $\m^\beta$ on $V\cap \t^*$ as simply
the image of $\n^\beta$ under the isomorphism $R_\g$:
$$ \m^\beta=R_\g(\n^\beta)\in D'(V\cap\t^*)^{W-\on{alt}}.$$
A more geometric construction of $\m^\beta$ is obtained as follows. 
Consider the Hamiltonian $T$-space $(M,\om,\Phi_T)$ where the $T$-moment map 
$\Phi_T$ is $\Phi$ followed by projection $\g^*\to \t^*$. 
Then $\m^\beta$ is the DH-distribution for this 
Hamiltonian $T$-space, corresponding to the $T$-equivariant 
cocycle $(\tpi)^{n_+} \tau_M \kappa_T(\beta)$ where $\tau_M=\Phi^*\tau$:  
\begin{equation}
\label{eq:defm0b}
 \m^\beta=(\tpi)^{n_+} \sum_I \big(\f{1}{\tpi}\f{\p}{\p \mu}\big)^I
(\Phi_{T})_*\big(\tau_{M} \kappa_T(\beta) e^{\tpi \om}\big)^{[\dim M]}_I.
\end{equation}
If $\beta$ is compactly supported, this can be written in terms of Fourier 
transform:  
$$ \l \mf{m}^\beta,e^{-\tpi\l\cdot,\xi\r}\r= (\tpi)^{n_+}\int_M \beta(\xi)
\tau_M(\xi) e^{\tpi (\om-\l\Phi_T,\xi\r)},\ \ \xi\in\t.$$
Due to the factor $\tau_M$, the integral is localized to an arbitrary 
small neighborhood of $X=\Phinv(V\cap \t^*)$. That is, it depends only 
on the restriction of $\beta$ to an arbitrarily small neighborhood 
of $X$. If the moment map $\Phi$ is transverse to $\t^*$, so that 
$X$ is a smooth submanifold, this becomes more concrete: 
\begin{proposition}\label{prop:transverse}
Suppose the moment map $\Phi$ is transversal to $\t^*$, 
and let $X=\Phinv(V\cap \t^*)$  with 
moment map $\Phi_X=\iota^*\Phi$ and 2-form $\om_X=\iota^*\om$ 
given by pull-back under the inclusion $\iota:\,X\hra M$. 
Then the distribution 
$\m^\beta$ for a cocycle $\beta\in\Om_G(M)$ is 
just the DH-distribution for the Hamiltonian $T$-space
$(X,\om_X,\Phi_X)$, associated to the cocycle
$\iota_X^*\beta$. 
\end{proposition}
\begin{proof}
In the transverse case, a tubular neighborhood of $X$ in $M$ may be 
identified in a $T$-equivariant way with $X\times\mf{p}$, in such a way 
that the diagram 
$$
\xymatrix{ 
X\times\mf{p} \ar[r]\ar[d]&M\ar[d]_\Phi\\ 
(V\cap T)\times \mf{p}\ar[r] & V}
$$
commutes. Then $\tau_M=\Phi^*\tau $ represents the Thom class of 
the vector bundle $\pi_X:\, X\times\mf{p}\to X$. The distribution 
$\m^\beta$ does not change if  
$\kappa_T(\beta) e^{\tpi \om}$ is replaced by 
$\pi_X^*(\beta_X e^{\tpi \om_X})$, since the difference is 
$\d_T$-exact. But then, using $\Phi_T=\Phi_X\circ \pi_X$,  
$$ (\Phi_T)_*(\tau_M \kappa_T(\beta) e^{\tpi \om})
=(\Phi_X)_* (\pi_X)_*(\tau_M \pi_X^*(\beta_X e^{\tpi \om_X})
=(\Phi_X)_* (\beta_X e^{\tpi \om_X}).$$
\end{proof}

We refer to Woodward \cite{wo:cl} for interesting examples of
Hamiltonian $G$-spaces with $\Phi$ transverse to $\t^*$. In the case
of a regular coadjoint orbit, Proposition \ref{prop:transverse} is
illustrated by Equation \eqref{eq:coadjoint}.

\section{DH-distributions for q-Hamiltonian $G$-spaces}
\label{sec:qdh}
Throughout this Section, we assume that $\rho$ (the half-sum of
positive roots) is a weight of $G$.  For instance, $G$ may be a
product of a simply connected group and a torus, or $G=\U(n)$ with $n$
odd. We identify $\g\cong\g^*$ by means of an invariant inner product
on $\g$. Similar to the Hamiltonian setting, any equivariant cocycle
$\beta\in\Om_G(M)$ on a q-Hamiltonian $G$-space $(M,\om,\Phi)$ defines
an equivariant current $\Phi_*(\beta e^{\tpi \om})$ on $G$, closed
under $\d_G'=\d_G-\tpi \eta_G$. Our first goal is, therefore, to
associate to any $\d_G'$-cocycle on $G$ an invariant distribution.

\subsection{The isomorphism $R_G$}
\label{subsec:rg}
The group analogue to the product of positive roots 
$\mu\mapsto \prod_{\alpha\in\mf{R}_+} 2\pi \alpha\cdot\mu$ 
is the function $t\mapsto i^{-n_+}A(t)$, 
where 
$$ A=\sum_{w\in W}(-1)^{l(w)} \eps_{w\rho}$$ 
is the Weyl denominator. Indeed, 
$$  i^{-n_+} A(\exp\mu)=\prod_{\alpha\in\mf{R}_+}2\sin(\pi \alpha\cdot\mu).$$
Pick orientations on $G,T$, compatible with the given orientation on
$\mf{p}$, and let $\d\vol_G$ and $\d\vol_T$ be the volume forms
defined by the orientation and the inner product. 
Let $V\subset G$ be an invariant open subset. Similar to Lemma
\ref{lem:restriction} we have:
\begin{lemma}\label{lem:rg}
There is a unique linear isomorphism $R_G:\,\D'(V)^G\to 
\D'(V\cap T)^{W-\on{alt}}$, with the following two properties: 
\begin{enumerate}
\item \label{it:lin1} {\bf ($C^\infty$-linearity)}
$R_G(f\mf{n})=f|_{V\cap T}R_G(\n)$ for all $f\in C^\infty(V)^G$
and $\n\in \D'(V)^G$. 
\item\label{it:ft1} {\bf (Fourier coefficients)}
For compactly supported $\n$, the Fourier
transform of $R_\g\n$ is given by
$$
\l R_G(\n),\,\ol{\eps_{\lambda+\rho}}\r =\f{i^{-n_+}}{\vol_{G/T}}\ \l
\n,\ol{\chi_\lambda}\r,\ \lambda\in \Lambda^*_+.$$
\end{enumerate}
For the volume form on $G$ one finds, 
$$ R_G(\d\vol_G)=i^{-n_+}A\ \d\vol_T$$
\end{lemma}
\begin{proof}
Define $R_G$ as the inverse dual map to the isomorphism,
$$ C^\infty(V)^G_{\on{comp}}\to C^\infty(V\cap
T)^{W-\on{alt}}_{\on{comp}} ,\ f\mapsto \f{\vol_{G/T}}{\# W}\ i^{-n_+}
A\, f|_{V\cap T} $$
Then $R_G$ clearly satisfies \eqref{it:lin1}, while Equation \eqref{it:ft1}
follows from the Weyl character formula:
\beq \l\n,\ol{\chi_\lambda}\r&=&\f{\vol_{G/T}}{\# W}\ 
i^{n_+} \l R_G\n,\,\ol{ A\,\chi_\lambda|_T}\r \\
&=&\f{\vol_{G/T}}{\# W}\ i^{n_+} \sum_{w\in W} (-1)^{l(w)} \l R_G\n,
\ol{\eps_{w(\lambda+\rho)}}\r\\
&=& \vol_{G/T}\,i^{n_+} \l R_G\n,\ol{\eps_{\lambda+\rho}}\r.\eeq
Similarly, the formula for
$R_G(\d\vol_G)$ is a consequence of the Weyl integration formula.
\end{proof}
The maps $R_\g$ and $R_G$ are related as follows: 
\begin{lemma}\label{lem:duflo}
The Duflo map $\on{Duf}_G:\,\D'(\g)^G_{\on{comp}}
\to \D'(G)^G,\ \n\mapsto \exp_* J^{1/2}\n$ 
intertwines the two maps $R_\g$ and $R_G$. That is, 
$$ R_G\circ \on{Duf}_G= \exp_*\circ R_\g$$
\end{lemma}
\begin{proof}
The identity $
i^{-n_+} A(\exp\mu)= J^{1/2}(\mu)\ \prod_{\alpha\in\mf{R}_+}
2\pi \alpha\cdot\mu$
for $\mu\in\t$ shows that for any test function $f\in C^\infty(G)^G$, 
$$ \exp^* (i^{-n_+}A\,f|_T)=
\, (\prod_{\alpha\in\mf{R}_+}2\pi \alpha)\ J^{1/2}\exp^* f.$$
Using the definition of $R_\g,R_G$ as dual maps to maps on test 
functions, this yields $(R_G)^{-1}\circ \exp_*
=\Duf_G \circ (R_\g)^{-1}$. 
\end{proof}
We will now use $R_G$ to construct a map $\ca{C}_G(V)\to \D'(V)^G$,
vanishing on cocycles for $\d_G'=\d_G-\tpi \eta_G$.  This map will 
be given in terms of a commutative diagram 
\begin{equation}\label{eq:diagram}
 \xymatrix {\ca{C}_G(V)\ar[r]\ar[d] &\D'(V)^G\ar[d]_{R_G}^\cong\\
\ca{C}_T(V\cap T)^{W-\on{alt}}\ar[r] & \D'(V\cap T)^{W-\on{alt}}}
\end{equation}
Here the $W$-action on $\ca{C}_T(V\cap T)$ is $(w.\phi)(\xi)=
w_*(\phi(w^{-1}.\xi))$. The lower horizontal map in the diagram 
is given by 
\begin{equation}\label{eq:hhh}
 \ca{C}_T(V\cap T)\to \D'(V\cap T),\ \ 
\psi \mapsto \m=\sum_I 
(\f{1}{\tpi}\f{\p}{\p \mu})^I \psi_I^{[\dim \t]}
\end{equation}
where we use $(\f{1}{\tpi}\f{\p}{\p \mu})^I$ to denote both a constant
coefficient differential operator on $\t$, and also the induced
differential operator on $T$. Note that this map is $W$-equivariant
and vanishes on $\d_T'=\d-\tpi\eta_T$-cocycles.
If $\psi$ has compact support, this map can be characterized 
in terms of Fourier coefficients by  
\begin{equation}\label{eq:fcoef}
 \l\m,\ol{\eps_\lambda}\r=\l\psi(\lambda),\ol{\eps_\lambda}\r,\ \ \lambda\in\Lambda^*.\end{equation}

To construct the left vertical map, consider the $T$-equivariant
tubular neighborhood $\pi:\,U\to T$ of $T\subset G$, described in
Section \ref{subsec:abel}.  Recall the 2-form $\gamma\in \Om^2(U)^T$,
satisfying $\d_T\gamma=\kappa_T(\eta_G)-\pi^*\eta_T$.  Restrict the bundle
$\pi:\,U\to T$ to the intersection $V\cap T$, and let $\tau\in
\Om_T^{2n_+}(\pi^{-1}(V\cap T))$ be a representative of the
$T$-equivariant Thom class, supported in the intersection
$\pi^{-1}(V\cap T)\cap V$. 
Then
\begin{equation}\label{eq:realize}
 \ca{C}_G(V)\to \ca{C}_T(V\cap T),\ \ 
\phi\mapsto \psi=(\tpi)^{n_+} \pi_*(\tau\,e^{\tpi \gamma}\,\kappa_T(\phi))
\end{equation}
is a well-defined map intertwining $\d_G'=\d_G-\tpi\eta_G$ with
$\d_T'=\d_T-\tpi\eta_T$.  Since the action of $w\in W$ on $\mf{p}$ changes
the orientation by $(-1)^{l(w)}$, we may choose $\tau$ in such a way
that for $g\in N_G(T)$, $g^*\tau(\xi)=(-1)^{l(w)}\tau(w^{-1}.\xi)$.
With this choice, the map \eqref{eq:realize} takes values in 
$\ca{C}_T(V\cap T)^{W-\on{alt}}$. We define the map 
$\ca{C}_G(V)\to \D'(V)^G,\ \phi\mapsto \n$ in the unique way making
\eqref{eq:diagram} commute. That is,
 \begin{equation}\label{eq:composed}
\n= (\tpi)^{n_+} R_G^{-1}\Big(\sum_I
(\f{1}{\tpi}\f{\p}{\p \mu})^I
\pi_*(\tau\,e^{\tpi \gamma}\,\kappa_T(\phi))_I^{[\dim \g]}\Big).
\end{equation} 
By construction, \eqref{eq:composed} vanishes on $\d_G'$-coboundaries, and the
induced map in $\d_G'$-cohomology is independent of the choice of
$\tau$.

\subsection{DH-distributions for q-Hamiltonian G-spaces}
Let $(M,\om,\Phi)$ be an oriented q-Hamiltonian $G$-space,
possibly degenerate, with proper moment map 
$\Phi:\,M\to V\subset G$.
 (Recall that if $\om$ is minimally degenerate, 
there exists a distinguished orientation, see Lemma \ref{lem:ori}).
Let $\beta\in \Om_G(M)$ be a closed 
equivariant differential form.  Then the current 
$$ \phi=\Phi_*(\beta e^{\tpi \om})\in C_G(V)$$
is $\d_G'=\d_G-\tpi \eta_G$-closed. Let $\n^\beta\in\D'(V)^G$ be
defined by \eqref{eq:composed}.  The distribution
$\m^\beta=R_G(\n^\beta) \in \D'(V\cap T)^{W-\on{alt}}$ admits an
interpretation similar to the Hamiltonian case, Section
\ref{subsec:mbeta}. Indeed, let $(N,\Om_N,\Phi_N)$ be the
q-Hamiltonian $T$-space given as the `Abelianization' of
$(M,\om,\Phi)$. Then $\m^\beta$ is the DH-distribution corresponding
to the cocycle $(\tpi)^{n_+} \tau_N \kappa_T(\beta)\in \Om_T(N)$ where
$\tau_N=\Phi^*\tau$:
\begin{equation}\label{eq:push}
 \m^\beta=(\tpi)^{n_+}\,\sum_I \big(\f{1}{\tpi}\f{\p}{\p \mu}\big)^I
(\Phi_N)_*\Big(\tau_N\kappa_T(\beta) e^{\tpi \om_N}\Big)_I^{[\dim M]}.
\end{equation}
If $\Phi$ is transverse to the maximal torus, the same argument 
as for Proposition \ref{prop:transverse} shows that $\m^\beta$ 
is a DH-distribution for the q-Hamiltonian $T$-space
$(X,\om_X,\Phi_X)$ where $X=\Phinv(T)$ and $\om_X,\Phi_X$ are pull-backs 
of $\om,\Phi$.

\begin{theorem}[Properties of the q-Hamiltonian DH-distributions]
\label{th:mb}
Let $\beta\in \Om_G(M)$ be an equivariant cocycle, and $\n^\beta$ 
the associated distribution. 
\begin{enumerate}
\item\label{it:indep}
$\n^\beta$ depends only on the cohomology class of $\beta$.
\item\label{it:loc}
If $V'\subset V$ is an
invariant open subset, then $\n^\beta|_{V'}$ is
the DH-distribution corresponding to 
$\beta|_{\Phinv(V')}$.  
\item\label{it:smooth}
If $g\in G$ is a $\beta$-regular value of
$\Phi$, the distribution $\n^\beta$ is smooth at $g$.
\item \label{it:integrand} 
If $\beta$ has compact support, then the
distribution $\n^\beta$ has compact support contained in
$\Phi(\on{supp}(\beta))$.  Its Fourier coefficients
$\l\n^\beta,\ol{\chi_\lambda}\r$ 
are given by the integral 
$$
\l\n^\beta,\ol{\chi_\lambda}\r=(-2\pi)^{n_+} \vol_{G/T}\,\int_N \tau_N(\lambda+\rho) 
\beta(\lambda+\rho)e^{\tpi \om_N}
\Phi_N^*\ol{\eps_{\lambda+\rho}},\ \ \lambda\in\Lambda^*_+ 
$$
which localizes to the fixed point set of $(\lambda+\rho)_M$:
$$
\l\n^\beta,\ol{\chi_\lambda}\r=\dim V_\lambda\, 
\sum_{F\in \ca{F}(\lambda+\rho)}\int_F \f{\iota_F^* \beta(\lambda+\rho)\,e^{\tpi \om_F}
\Phi_F^*\ol{\eps_{\lambda+\rho}}}{\Eul(\nu_F,\lambda+\rho)}.
$$
Here $\om_F\in \Om^2(F)$ and $\Phi_F:\,F\to T$ are defined as
pull-backs of $\om$ and $\Phi$ to the fixed point manifold $F$.
\end{enumerate}
\end{theorem}
\begin{proof}
 \eqref{it:indep}, \eqref{it:loc} and \eqref{it:smooth} follow from
the properties of the map $\ca{C}_G(V)\to \D'(V)^G$. In particular, if
$g\in G$ is a $\beta$-regular value, the current $\Phi_*(\beta e^{\tpi
\om})$ is smooth near $g$, hence so is the distribution $\n^\beta$.

Suppose next that $\beta$ is compactly supported.  Using
$(\f{1}{\tpi}\f{\p}{\p\mu})^I\eps_\lambda=\lambda^I\eps_\lambda$, the
Fourier coefficients of the distribution $\m^\beta$ may 
be written, 
$$
\l\m^\beta,\ol{\eps_{\lambda+\rho}}\r=(\tpi)^{n_+} \int_N \tau_N(\lambda+\rho) 
\beta(\lambda+\rho)e^{\tpi \om_N}
\Phi_N^*\ol{\eps_{\lambda+\rho}},\ \ \lambda\in\Lambda^*_+. 
$$
This implies the first Formula in \eqref{it:integral}, by Lemma
\ref{lem:rg}\eqref{it:ft1}.  Since the integrand is
$\d_{\lambda+\rho}$-closed, the integral may be computed by
localization, see 
Theorem \ref{th:bv}. The fixed
point components $F$ are all contained in $\Phi^{-1}(T)$, by
equivariance of the moment map. The second Formula in \eqref{it:integral}
follows since $\iota_F^*\Phi_N=\iota_F^*\Phi$ and $\iota_F^*(\om_N)=\iota_F^*\om$, 
and since $\iota_T^*\tau$ is cohomologous
to
$$(-1)^{n_+}\prod_{\alpha\in\mf{R}_+}\alpha\cdot(\lambda+\rho)
=(-2\pi)^{-n_+}\f{\dim V_\lambda}{ \vol_{G/T}},$$
where we used the Weyl dimension formula.
\end{proof}

\subsection{Relation to intersection pairings}\label{sec:intersect}
Our goal in this Section is to prove: 
\begin{theorem}\label{th:interpret}
Let $V\subset G$ be an invariant open subset containing the group unit 
$e$, and $(M,\om,\Phi)$ a q-Hamiltonian $G$-space, with proper moment
map $\Phi:\,M\to V\subset G$. Suppose $\beta\in
\Om_G(M)$ is an equivariant cocycle, and that $e\in V$ is a
$\beta$-regular value of $\Phi$. Then the distribution $\n^\beta$ is
smooth near $e$, and
$$ \n^\beta(e)=(\tpi)^{\dim G}\vol_{G}\sum_j \f{1}{k_j} 
\int_{(M\qu G)_j} \beta_{\on{red}}e^{\tpi \om_{\on{red}}}.
$$
Here the sum is over connected components of $M\qu G$, and $k_j$ is
the cardinality of a generic stabilizer for the $G$-action on the
$j$th component of $\Phinv(e)$. 
\end{theorem}
The idea of proof is to compare with the DH distribution for the
Hamiltonian $G$-space described in \ref{subsec:log}.  We first discuss
a similar problem for equivariant currents. Recall that 
in \eqref{eq:varpi}, we defined a  
2-form $\varpi\in\Om^2(\g)^G$ with $\d_G\varpi=\exp^*\eta_G-\eta_\g$.  

\begin{lemma}\label{lem:smooth1}
Suppose $V_0\subset \g$ is an invariant open subset of $\g$, such that
$\exp$ restricts to a diffeomorphism onto the image $V=\exp(V_0)$.
Let $\phi\in\Om_G(V)\subset \ca{C}_G(V)$ be a smooth equivariant
current, and $\phi_0:=e^{-\tpi\varpi} \exp^*\phi\in \Om_G(V_0)$.  Then
$\phi_0$ is closed under $\d_G-\tpi\eta_\g$ if and only if $\phi$ is
closed under $\d_G-\tpi\eta_G$. The corresponding (smooth) measures
$\n_0$ on $V_0$ and $\n$ on $V$ are related by the Duflo map:
\begin{equation}\label{eq:measures}
 \n=\exp_*(J^{1/2}\n_0).
\end{equation}
\end{lemma}
\begin{proof}
The first claim is clear since conjugation by
$e^{\tpi\varpi}$ intertwines $\d_G-\tpi\exp^*\eta_G$ with
$\d_G-\tpi\eta_\g$.  
For the second claim, it suffices, by Corollary \ref{cor:cutoff}, 
to consider the case that $\phi,\phi_0$ have compact support. 
Let $\m=R_G(\n)$ and $\m_0=R_\g(\n_0)$. By Lemma \ref{lem:duflo} 
we have to show 
$$\m=\exp_*\m_0.$$
Since $\m_0,\m$ are compactly supported and $W$-alternating, 
it is enough to show that 
\begin{equation}\label{eq:fcoeff}
 \l\m,\ol{\eps_{\lambda+\rho}}\r=\l \m_0,e^{-\tpi
 \l\cdot,\lambda+\rho\r}\r,\ \ \lambda\in\Lambda^*_+.
\end{equation}
By definition
$$ \l\m,\ol{\eps_{\lambda+\rho}}\r=(\tpi)^{n_+}\int_{U}
\tau (\lambda+\rho) 
e^{\tpi \gamma} \phi(\lambda+\rho) \pi^*\ol{\eps_{\lambda+\rho}}.$$
where the integrand is closed for the differential
$\d_{\lambda+\rho}$. Integrating over the fibers of $\pi:\,U\to T$, 
and using $\iota_T^*\gamma=0$,
this gives
$$ \l\m,\ol{\eps_{\lambda+\rho}}\r=(\tpi)^{n_+}
\int_T \iota_T^* \phi'(\lambda+\rho)\ol{\eps_{\lambda+\rho}}.$$
A similar argument for 
$\phi_0$ shows that  
$$ \l\m_0,e^{-\tpi \l\cdot,\xi\r}\r=
(\tpi)^{n_+} \int_\t \iota_\t^*\phi_0(\xi)e^{-\tpi \l\cdot,\xi\r},\ \ 
\xi\in\on{int}(\t_+).$$
But $\exp_T^* \iota_T^*\phi=\iota_\t^*\phi_0$, 
since $\iota_\t^*\varpi=0$. Hence we have proved \eqref{eq:fcoeff} and 
therefore the Lemma. 
\end{proof}

\begin{proof}[Proof of Theorem \ref{th:interpret}]
Replacing $V$ with a smaller neighborhood of $e$ if necessary (and $M$
with the corresponding pre-image under $\Phi$), we may assume that all
points of $V$ are $\beta$-regular values of $\Phi$. Furthermore, we
may assume that $V$ is the diffeomorphic image under $\exp$ of an
invariant open neighborhood $V_0\subset \g$. Let $(M,\om_0,\Phi_0)$ be
given by linearization (cf. Section \ref{subsec:lin}). By
construction, the currents $\phi=\Phi_*(\beta e^{\tpi \om})$ on $V$
and $\phi_0=(\Phi_0)_*(\beta e^{\tpi \om_0})$ are smooth, and
$\phi=\exp_*(e^{\tpi \varpi} \phi_0)$. Hence, Lemma \ref{lem:smooth1}
shows that the measures $\n_0^\beta$ and $\n^\beta$ are related by the
Duflo map $\on{Duf}_G:\,\D'(V_0)^G\to \D'(V)^G$. Since $J^{1/2}(0)=1$,
it follows in particular that
$$ \n^\beta(e)=\n_0^{\beta}(0)$$
where we use the Riemannian measures on $\g$ and on $G$ to identify
smooth measures with functions. Since
$\om_{\on{red}}=(\om_0)_{\on{red}}$, the Theorem follows from 
the interpretation of DH-distributions for Hamiltonian spaces, Proposition
\ref{prop:nb0}\eqref{it:interpret}.
\end{proof}

\begin{remark}
More generally, the value of $\n^\beta$ at any $\beta$-regular value
$g$ of $\Phi$ is given by an integral over components of $\Phinv(g)/G_g$, 
similar to \ref{prop:nb0}\eqref{it:interpret}, with the symplectic volume of a
(co-)adjoint orbit $G.\mu$ replaced by the q-Hamiltonian volume of a
conjugacy class \eqref{eq:covol}. For $g=\exp\mu$ with $\mu$ sufficiently
small, this follows by the same argument as for $g=e$, using that
$\Vol(G.\exp\mu)=J^{1/2}(\mu) \Vol(G.\mu)$. For general $g$, the
result can be obtained using cross-sections \cite{al:mom}.
\end{remark}

\subsection{Fixed point formula for intersection pairings}
Suppose $(M,\om,\Phi)$ is an oriented q-Hamiltonian $G$-space, with
proper moment map $\Phi:\,M\to G$. (In particular, $M$ is compact.)
Assume that $e\in G$ is a regular value of $\Phi$. 

Given an equivariant cocycle $\beta\in \Om_G(M)$ consider the Fourier
expansion of $\n^\beta$:
$$ \n^\beta=\f{1}{\vol_G}\,\sum_{\lambda\in\Lambda^*_+}
\l\n^\beta,\ol{\chi_\lambda}\r\,
\chi_\lambda$$
where $\d\vol_G$ is used to identify distributional measures on $G$ 
with distributional functions.  
Combining Theorem \ref{th:interpret} with
\ref{th:mb}\eqref{it:integrand}, and using $\chi_\lambda(e)=\dim
V_\lambda$, one recovers the following result from \cite{al:gr}:
\beq
\lefteqn{
\sum_j \f{1}{k_j} 
\int_{(M\qu G)_j} \beta_{\on{red}}e^{\tpi \om_{\on{red}}}}
\\
&=&(\tpi)^{-\dim G}
 \sum_{\lambda\in\Lambda^*_+}
(\f{\dim V_\lambda}{\vol_G})^2 \sum_{F\in\F(\lambda+\rho)}
\int_F \f{\iota_F^* \beta(\lambda+\rho)\,e^{\tpi \om_F}
}{\Eul(\nu_F,\lambda+\rho)}\Phi_F^*\ol{\eps_{\lambda+\rho}}
\eeq
Of course, the formal substitution $g=e$ is only valid if 
the Fourier coefficients decay sufficiently fast, to ensure 
convergence of this series. In the general case, one can introduce 
a convergence factor $\exp(-\eps||\lambda+\rho||^2)$ in the sum, and
obtains an equality for  $\eps\to 0$. See \cite{li:he} and \cite{al:du}
for more detailed discussion. 

\begin{remark}
If $G$ is simply connected, and $M$ a compact, connected
q-Hamiltonian $G$-space with a minimally degenerate 2-form $\om$, then
the fibers of the moment map $\Phi$ are connected \cite{al:mom}.  In
this case, the formula directly gives intersection pairings on $M\qu
G$ in terms of fixed point contributions.
\end{remark}

\section{DH-distributions for higher q-Hamiltonian forms}
\subsection{Currents associated to higher q-Hamiltonian forms}
Suppose $(M,\om,\Phi)$ is an oriented q-Hamiltonian $G$-space, 
with proper moment map $\Phi:\,M\to V\subset G$, and let 
$\beta\in \Om_G(M)$ be an equivariant cocycle.  
Given an invariant polynomial $p\in \on{Pol}(\g)^G$, suppose 
$\om^p\in\Om_G(M)$ is a higher q-Hamiltonian form, i.e. satisfying 
$\d_G\om^p=\Phi^*\eta_G^p$. Then 
\begin{equation}\label{eq:pphi}
\phi=\Phi_*(\beta e^{\tpi \om^p})
\end{equation}
is closed with respect to the differential $\d_G-\tpi\eta^p_G$.
Note however that $\beta(\xi) e^{\tpi \om^p(\xi)}$ does not lie in 
$\Om_G(M)$, in general, since the exponential does not depend 
polynomially on $\xi\in \g$. On the other hand, we have to insist 
on polynomial dependence in order for formulas such as 
\eqref{eq:composed} to make sense. 

To get around this difficulty, we 
take $p$ to be an invariant polynomial of the form
\begin{equation}\label{eq:pform}
p(\xi)=p_1(\xi)+\sum_{j=2}^l \delta_j p_j(\xi)
\end{equation}
where $p_1(\xi)=\hh ||\xi||^2$, the $p_j$ for $j\ge 2$ are invariant
polynomials, and $\delta=(\delta_2,\ldots,\delta_l)$ are formal
parameters. That is, instead of a general invariant polynomial $p$ we 
consider only deformations of the quadratic polynomial. 
It will be convenient to introduce a notation for the perturbation term, 
$$ q(\xi)=\sum_{j=2}^l \delta_j p_j(\xi).$$ 
We also assume that the leading term of $\om^p$ is the given 2-form
$\om$, so that
\begin{equation}\label{eq:omp}
 \om^p=\om+\om^q=\om+\sum_{j=2}^l \delta_j \om^{p_j}
\end{equation}
where $\om^{p_j}\in\Om_G(M)$ are higher q-Hamiltonian forms
corresponding to $p_j$. Then $\beta e^{\tpi \om^p}$ is defined as an
element of $\Om_G(M)[[\delta]]$ where
$\C[[\delta]]=\C[[\delta_2,\ldots,\delta_l]]$ denotes the ring of
formal power series. Similarly, Equation \eqref{eq:pphi}
defines an element of $\ca{C}_G(V)[[\delta]]$, closed under the 
differential $\d_G-\tpi\eta^p_G$.

\subsection{Witten's change of variables}
Suppose $\phi\in \ca{C}_G(V)[[\delta]]$ is closed under 
the differential $\d_G-\tpi\eta^p_G$. We would like to associate to 
$\phi$ an invariant distribution $\n\in \D'(V)^G[[\delta]]$. 

As a first step, we define a current on $V\cap T$, 
by a formula similar to \eqref{eq:realize}
\begin{equation}\label{eq:current}
\psi=(\tpi)^{n_+} \pi_*(\tau e^{-\tpi\gamma^p}\, \kappa_T(\phi))
\in \ca{C}_T(V\cap T)[[\delta]], 
\end{equation}
where $\gamma^p$ is given by \eqref{eq:gammap}.  Note that $\psi$ is
well-defined since
$e^{-\tpi \gamma^p}=e^{-\tpi \gamma}e^{-\tpi \gamma^q}$ lies in
$\Om_T(U)[[\delta]]$, i.e. the coefficients depend {\em polynomially}
on the Lie algebra variables. Furthermore, $\psi$ is closed for the
differential $\d-\tpi \eta_T^p$. 

The second step associates to any $\d-\tpi \eta_T^p$-cohomology class
a distribution on $V\cap T$. For this, we cannot directly use the map
\eqref{eq:hhh}, since this map does not vanish on $\d-\tpi
\eta_T^p$-coboundaries unless $p=\hh ||\xi||^2$. The underlying
problem is that the current $\psi(\lambda)\ol{\eps_\lambda}$ (for weights
$\lambda \in\Lambda^*$) is not $\d$-closed, in general. Instead, using
$(\d+\tpi \lambda\cdot \theta_T)\ol{\eps_\lambda}=0$ and $(\d-\tpi
p'(\xi)\cdot \theta_T)\psi(\xi)=0$, we have:
\begin{lemma}
For $\xi\in\t$ and $\lambda\in \Lambda^*$, the current
$\psi(\xi)\ol{\eps_\lambda}\in \ca{C}(V\cap T)$ is $\d$-closed if and
only if $p'(\xi)=\lambda$.
\end{lemma}
This observation motivates a `change of variables'
$\xi=(p')^{-1}(\lambda)$, as in Witten's paper \cite[Equation
(5.9)]{wi:tw}.  Let $\on{Pol}(\g,\g)\subset C^\infty(\g,\g)$ be the
algebra of polynomial maps from $\g$ to itself. The transformation 
$\xi\mapsto p'(\xi)$ 
is a well-defined
invertible element of the algebra $\on{Pol}(\g,\g)[[\delta]]$, with
leading term the identity map. (See Appendix \ref{app:fcov} for more
on the change of variables). We will need the following fact regarding
the Jacobian $p''$ of the change of variables. View $p''$ as an
element of $\on{Pol}(\g,\,\on{End}(\g))[[\delta]]$, given in terms of
an orthonormal basis $e_a$ of $\g$ and the associated coordinates
$\xi^a$ by the matrix $\f{\p^2 p}{\p \xi^a \p \xi^b}$.  
\begin{lemma}\label{lem:alemma}
For any $\xi\in\t$, the linear map $p''(\xi)\in \on{End}(\g)[[\delta]]$
preserves the decomposition $\g=\t\oplus \mf{p}$. Furthermore, 
if $\xi$ is a regular element of $\t$ 
the determinant of $p''(\xi)|_{\mf{p}}$ is given by 
$$ 
(\det\, p''(\xi)|_{\mf{p}})^{1/2}=
\prod_{\alpha\in\mf{R_+}}
\f{\alpha\cdot p'(\xi)}{\alpha\cdot \xi}.$$
\end{lemma}
\begin{proof}
Let us treat the $\delta_i$ as real variables rather than as formal 
parameters. For sufficiently small $\delta_i$, the derivative $p':\g\to\g$ is a
well-defined diffeomorphism of neighborhoods of $0$
containing $\xi$ resp. $p'(\xi)$.  By invariance, the derivatives of
$p$ at $\xi$ vanish in ${\mf{p}}$-directions, which implies that
$p''(\xi)$ has block form
$$p''(\xi)=\left(\begin{array}{cc}p''(\xi)|_{\t}
&0\\0 & p''(\xi)|_{{\mf{p}}}\end{array}\right).$$
Suppose now that $\xi$ is a regular element. 
The matrix $p''(\xi)|_{{\mf{p}}}$ is the Jacobian for 
the transformation 
$$T_\xi(G\cdot\xi)\cong {\mf{p}}\to T_{p'(\xi)}(G\cdot p'(\xi))\cong
{\mf{p}}
$$
induced by $p'$. Hence its determinant is the ratio of the Riemannian volumes 
of the orbits through $\xi$ resp. $p'(\xi)$, with the Riemannian 
metric induced from $\g$. The Lemma follows, since 
the Riemannian volume 
of an adjoint orbit $G\cdot\zeta\subset \g$ through a regular element
$\zeta\in\t$  is equal to $|\det_{\mf{p}}(\ad_\zeta)|\,\Vol_{G/T}$
where 
$$ |{\det}_{\mf{p}}(\ad_\zeta)|=(\prod_{\alpha\in\mf{R}_+}2\pi \alpha\cdot \zeta)^2.$$
\end{proof}

Let $\mf{A}_\g$ denote the algebra automorphism of
$\on{Pol}(\g)[[\delta]]$ given by
$$ (\mf{A}_\g F)(p'(\xi))=F(\xi).$$
Since the change of variables operator 
$\mf{A}_\g$ commutes with the adjoint action, it induces 
algebra automorphisms of the Cartan complexes $\Om_G(M)[[\delta]]$ 
and $\ca{C}_G(M)[[\delta]]$, for any $G$-manifold $M$. Let
$\mf{A}_\t$ be defined similarly, using the restriction $p|_\t$.

\begin{lemma}
The automorphism $\mf{A}_\t$ of $\ca{C}_T(V\cap T)[[\delta]]$
intertwines the differential $\d-\tpi \eta^p_T$ with $\d-\tpi
\eta_T$. In particular, if $\psi\in \ca{C}_T(V\cap T)[[\delta]]$ is
closed under $\d-\tpi \eta^p_T$, then $\mf{A}_\t \psi$ is closed under
$\d-\tpi \eta_T$.
\end{lemma}
\begin{proof}
It is obvious that $\mf{A}_\t$ commutes with $\d$ and intertwines
$\eta^p_T(\xi)=p'(\xi)\cdot\theta_T$ with $\eta_T(\xi)=\xi\cdot\theta_T$. 
\end{proof}
Hence, by composing $\mf{A}_\t$ with the map \eqref{eq:hhh}, we obtain
a map $\ca{C}_T(V\cap T)[[\delta]]\to \D'(V\cap T)[[\delta]]$,
vanishing on $\d-\tpi \eta^p_T$-coboundaries.  As it turns out, it
will be convenient to modify this map by the factor $S=(\det
p''|_{\mf{p}})^{1/2}$, and to
define a $W$-equivariant map $\ca{C}_T(V\cap T)[[\delta]]\to \D'(V\cap T)[[\delta]]$
by
\begin{equation} \label{hhh1}
\psi\mapsto \m=
\sum_I (\f{1}{\tpi}\f{\p}{\p \mu})^I (\mf{A}_\t (S \psi))_I^{[\dim \t]}.
\end{equation}
Hence there is a unique map $
\ca{C}_G(V)[[\delta]]\to \D'(V)^G[[\delta]]$ vanishing on $\d_G-\tpi
\eta^p_G$-coboundaries, for which the following diagram commutes:
\begin{equation}\label{eq:diagram2}
 \xymatrix{ \ca{C}_G(V)[[\delta]]\ar[r]\ar[d]&
\D'(V)^G[[\delta]]\ar[d]_{R_G}^\cong\\ \ca{C}_T(V\cap T)^{W-\on{alt}}[[\delta]] \ar[r]&
\D'(V\cap T)^{W-\on{alt}}[[\delta]]}
\end{equation}
Here the lower horizontal map is \eqref{hhh1} and the left vertical map 
is \eqref{eq:current}. 

\begin{lemma} Suppose $\phi\in\ca{C}_G(V)[[\delta]]$ is closed under 
$\d_G-\tpi \eta^p_G$, and let $\n\in \D'(V)^G[[\delta]]$ be its 
image under the above map. If $\phi$ is compactly supported, 
the Fourier coefficients of $\n$ are given by the formula, 
$$ \l\n,\ol{\chi_\lambda}\r=\f{\dim V_\lambda}
{(-1)^{n_+} \prod_{\alpha\in\mf{R}_+}\alpha\cdot\xi}
 \int_U
\tau(\xi)
e^{-\tpi\gamma^p(\xi)}\, \phi(\xi) \pi^*\ol{\eps_{\lambda+\rho}}
$$
where $\xi$ is the solution of $p'(\xi)=\lambda+\rho$. 
\end{lemma}
\begin{proof}
The Lemma follows from $\l\n,\ol{\chi_\lambda}\r=
i^{n_+} \vol_{G/T}\l\m,\ol{\eps_{\lambda+\rho}}\r$ and the 
calculation, 
\beq  \l\m,\ol{\eps_{\lambda+\rho}}\r&=&
S(\xi)\,\l \psi(\xi),\ol{\eps_{\lambda+\rho}}\r\\&=&
\prod_{\alpha\in\mf{R}_+}\f{\alpha\cdot(\lambda+\rho)}{\alpha\cdot\xi}
(\tpi)^{n_+} \l \tau(\xi)\,e^{-\tpi\gamma^p(\xi)}\,\phi(\xi),\pi^*\ol{\eps_{\lambda+\rho}}\r.
\eeq
using
$\prod_{\alpha\in\mf{R}_+}\,\alpha\cdot(\lambda+\rho)=(2\pi)^{-n_+}\f{\dim
V_\lambda}{\vol_{G/T}}$.
\end{proof}
Returning to the case of a q-Hamiltonian $G$-space, consider $\om^p$
of the form \eqref{eq:omp}, and an equivariant cocycle
$\beta\in\Om_G(M)$.  We define a generalized DH-distribution
$\n^\beta\in \D'(V)^G[[\delta]]$ as the image of the current
$\phi=\Phi_*(\beta e^{\tpi\om^p})\in \ca{C}_G(V)[[\delta]]$ under the map
\eqref{eq:diagram2}.

If $\beta$ is compactly supported (in particular if $M$ is compact), 
the Fourier coefficients of $\n^\beta$ can be expressed
in terms of the Abelianization $(N,\om_N,\Phi_N)$. 
Recall that by Proposition \ref{prop:abel}, the form 
$\om^p_N=\om^p-\Phi^*\gamma^p$ satisfies $\d_T \om^p_N=\Phi_N^*\eta_T^p$. 
We obtain: 
\begin{equation}
 \l\n^\beta,\ol{\chi_\lambda}\r=\f{\dim V_\lambda}
{(-1)^{n_+} \prod_{\alpha\in\mf{R}_+}\alpha\cdot\xi}
 \int_U
\tau(\xi)\beta(\xi) e^{\tpi\om^p_N(\xi)}\,  \pi^*\ol{\eps_{\lambda+\rho}}.
\end{equation}
Using localization, we find:
\begin{proposition}\label{prop:locform}
Suppose $\beta\in\Om_G(M)$ is a compactly supported cocycle. Then 
$$ \l\n^\beta,\ol{\chi_\lambda}\r=
\dim V_\lambda \sum_{F\in \F(\lambda+\rho)} \int_F
\f{\beta(\xi) e^{\tpi \om^p(\xi)}\Phi_F^*\ol{\eps_{\lambda+\rho}}}{\Eul(\nu_F,\xi)},\ \ \lambda\in\Lambda^*_+
$$
where $\xi$ is the solution of $p'(\xi)=\lambda+\rho$. 
\end{proposition}
\begin{proof}
We view the $\delta_j$ as real variables rather than just formal
parameters. For $\delta_j$ sufficiently small, the inverse
$\xi=(p')^{-1}(\lambda+\rho)$ is a well-defined smooth function of
$\delta_j$. Moreover, $\xi\to \lambda+\rho$ as $\delta_j\to 0$.  This
implies that $\xi_M^{-1}(0)\subset (\lambda+\rho)_M^{-1}(0)$, for
$\delta_j$ sufficiently small. Since the integrand in the definition of
$\m^\beta_{\lambda+\rho}$ is $\d_\xi$-closed, Proposition
\ref{prop:locform} follows from the localization
formula (Theorem \ref{th:bv}), using that the pull-back of
$\tau_M(\xi)$ to any $F$ is $T$-equivariantly cohomologous to
$(-1)^{n_+}\prod_{\alpha\in\mf{R}_+} \l\alpha,\xi\r$.
\end{proof}

\subsection{Interpretation}
By construction, the distribution $\n^\beta\in \D'(V)[[\delta]]$ is 
smooth at $\beta$-regular values of $\Phi$. Its restriction to the set of 
$\beta$-regular values $g\in V$ encodes intersection pairings on 
symplectic quotients $\Phinv(g)/G_g$. We will need the precise relationship 
only for $g=e$. Generalizing Theorem \ref{th:interpret} we will prove: 
\begin{theorem}\label{th:interpret1}
Let $(M,\om,\Phi)$ be a q-Hamiltonian $G$-space, with proper moment
map $\Phi:\,M\to V\subset G$, and $\beta\in\Om_G(M)$ is an equivariant
cocycle. Let $\om^p$ be a higher q-Hamiltonian form as in
\eqref{eq:omp}, and $\n^\beta\in \D'(V)[[\delta]]$ the DH-distribution
defined by these data.  Suppose $V$ contains the group unit $e$, and
is a $\beta$-regular value of $\Phi$. Then $\n^\beta$ is smooth near
$e$, and using $\d\vol_G$ to identify measures and functions,
$$ 
 \n^\beta(e)=(\tpi)^{\dim G}\vol_{G}\sum_j \f{1}{k_j} 
\int_{(M\qu G)_j} (\det(p'')\,\beta)_{\on{red}}e^{\tpi \om^p_{\on{red}}}.
$$
\end{theorem}
A combination of Theorem \ref{th:interpret1} and Proposition
\ref{prop:locform} (applied to $\n^{\beta_1}$ where
$\beta_1=\f{\beta}{\det(p'')}$) gives a localization formula for
intersection pairings:
\begin{theorem}[Localization formula for higher q-Hamiltonian forms]
\label{th:higher}
Let $(M,\om,\Phi)$ be a compact, connected q-Hamiltonian $G$-space, 
with group unit $e$ a regular value of the moment map. 
and $\beta\in \Om_G(M)$ an equivariant cocycle. Let
$p(\xi)=\hh ||\xi||^2+\sum_{j=2}^l \delta_j p_j(\xi)$ be an invariant 
polynomial and $\om^p=\om +\sum\delta_j \om^{p_j}\in\Om_G(M)$
a corresponding higher q-Hamiltonian form. Then 
\beq
\lefteqn{
\sum_j \f{1}{k_j} 
\int_{(M\qu G)_j} \beta_{\on{red}}e^{\tpi \om^p_{\on{red}}}}
\\
&=&\f{(\tpi)^{-\dim G}}{\vol_G^{2}}
 \sum_{\lambda\in\Lambda^*_+}
\f{(\dim V_\lambda)^2}{\det p''(\xi)}\sum_{F\in\F(\lambda+\rho)}
\int_F \f{\iota_F^* \beta(\xi)\,e^{\tpi \om_F^p(\xi)}
}{\Eul(\nu_F,\xi)}\Phi_F^*\ol{\eps_{\lambda+\rho}}
\eeq
where in each summand, $\xi\in\g[[\delta]]$ is given as the 
solution of $p'(\xi)=\lambda+\rho$. 
\end{theorem}

\subsection{Proof of Theorem \ref{th:interpret1}}
The idea of proof is to relate the distribution $\n^\beta$ to a
similar distribution for ordinary Hamiltonian spaces. Since we are
only considering a neighborhood of $e$, we may assume that $V$ is the
diffeomorphic image of an invariant neighborhood $V_0\subset \g$ of
$0$ under the exponential map. Furthermore, we may assume that $V$
consists of $\beta$-regular values of $\Phi$ so that the current
$\phi$ is in fact smooth.

We begin by considering arbitrary $\d_G-\tpi \eta^p_G$-closed currents
$\phi\in \ca{C}_G(V)[[\delta]]$, and the corresponding distributions
$\n\in\D'(V)[[\delta]]$ and $\m=R_G(\n)$. 
 
Let $\varpi^p\in\Om_G(\g)[[\delta]]$ be defined as in
\eqref{eq:varpip}.  Conjugation by $e^{\tpi \varpi^p}$ is an
automorphism of $\ca{C}_G(\g)[[\delta]]$, intertwining $\d_G-\tpi
\exp^*\eta^p_G$ with $\d_G-\tpi \eta^p_\g$. Hence the form
\begin{equation}
\label{eq:phi00} 
\phi_0=e^{-\tpi \varpi^p} \exp^*\phi\in
\ca{C}_G(V_0)[[\delta]]
\end{equation}
is closed for the differential $\d_G-\tpi \eta_\g^p$. The change of
variables operator $\mf{A}_\g$ intertwines
$\eta_\g^p(\xi)=-\d\l\cdot,p'(\xi)\r$ with $\eta_\g(\xi)=-\d\l\cdot,\xi\r$. 
Therefore $\mf{A}_\g\phi_0$ is closed for $\d_G-\tpi\eta_\g$, and 
we have the associated distribution, 
$$ \n_0=\sum_I \big(\f{1}{\tpi}\f{\p}{\p \mu}\big)^I(\mf{A}_\g
\phi_0)^{[\dim \g]}_I.$$
Similar to $\psi$, define a $\d_T-\tpi\eta_\t^p$-cocycle
\begin{equation}
\label{eq:psi00} 
\psi_0=(\tpi)^{n_+}(\pi_0)_*(\tau_0\,\kappa_T(\phi_0))\in\ca{C}_T(V_0\cap\t)[[\delta]],
\end{equation}
and the corresponding 
distribution, 
$$ \m_0=\sum_I \big(\f{1}{\tpi}\f{\p}{\p \mu}\big)^I (\mf{A}_\t
S\,\psi_0)^{[\dim \t]}_I.$$
\begin{lemma}
The distributions $\n_0$ and $\m_0$ are related by $\m_0=R_\g(\n_0)$. 
\end{lemma}
\begin{proof}
It suffices to prove this for compactly supported $\phi_0$. In this case, the 
Fourier transforms of $\n_0$ is $\l\n_0,e^{-\tpi\l\cdot,p'(\xi)\r}\r=
\l\phi_0(\xi),e^{-\tpi\l\cdot,p'(\xi)\r}\r$. For the Fourier transform of 
$\m_0$ we compute, 
\beq 
\l\m_0,e^{-\tpi\l\cdot,p'(\xi)\r}\r&=&
S(\xi) \l\psi_0(\xi),e^{-\tpi\l\cdot,\xi\r}\r\\
&=& (\tpi)^{n_+}\,S(\xi) \l \tau_0(\xi) \phi_0(\xi),e^{-\tpi\l\cdot,\xi\r}\r.\\
\eeq
Since $\phi_0$ is compactly supported, $\tau_0(\xi)$ may be replaced
by the cohomologous form of non-compact support 
$(-1)^{n_+}\prod_{\alpha\in\mf{R}_+}\l
\alpha,\xi\r$. This factor  combines
with the denominator of $S$, and we obtain 
$$
\l\m_0,e^{-\tpi\l\cdot,p'(\xi)\r}\r=(-1)^{n_+}\prod_{\alpha\in\mf{R}_+}\tpi
\l \alpha,p'(\xi)\r\ \l\n_0,e^{-\tpi\l\cdot,p'(\xi)\r}\r.$$
Now use the characterization of $R_\g$ by Fourier transforms, 
Lemma \ref{lem:restriction}\eqref{it:ft}.
\end{proof}

\begin{lemma}\label{lem:duflo2}
Suppose $\phi\in \Om_G(V)[[\delta]]$ is a smooth equivariant current,
closed under $\d_G-\tpi \eta^p_G$.  Let $\n$ and $\phi_0,\n_0$ be as
above. Then $\n,\n_0$ are related by the Duflo map:
$$\n=\exp_*(J^{1/2}\n_0).$$
\end{lemma}
\begin{proof}
Again, it is enough to prove this for compactly supported $\phi$.  We
have to show that $\m=R_G(\n)$ and $\m_0=R_\g(\n_0)$ are related by
$\m=(\exp_T)_*\m_0$.  As in the proof of Lemma \ref{lem:smooth1}, this
is verified by comparing the Fourier coefficients of the two measures:
First, defining $\xi$ by $p'(\xi)=\lambda+\rho$, the localization
formula (or integration over the fiber, using the properties of the Thom 
form) gives
$$\l\m,\ol{\eps_{\lambda+\rho}}\r=(\tpi)^{n_+} S(\xi)\, \int_U
\tau(\xi) \phi(\xi) e^{-\tpi\gamma^p(\xi)}
\pi^*\ol{\eps_{\lambda+\rho}}= S(\xi) \, \int_T \iota_T^*
\phi(\xi)\ol{\eps_{\lambda+\rho}} $$
where we have used that the integrand is $\d_\xi$-closed. The Fourier 
coefficient $\l\m_0',e^{-\tpi \l\cdot,\lambda+\rho\r}\r$ is computed 
similarly, and agrees with $\l\m',\ol{\eps_{\lambda+\rho}}\r$. 
\end{proof}

Above, we have seen that $\mf{A}\phi_0$ is a $\d_G-\tpi\eta_\g$-closed 
current. Another way of making $\phi_0$ closed under $\d_G-\tpi\eta_\g$
is to define 
$$\ti{\phi}_0(\xi)=e^{\tpi \l\cdot, q'(\xi)\r}\phi_0(\xi)\in 
\ca{C}_G(V_0)[[\delta]]. 
$$
Let $\ti{\n}_0\in \D'(V_0)[[\delta]]$ be the corresponding distribution. 
As it turns out, we will find it much easier to interpret $\ti{\n}_0$ 
rather than $\n_0$, in the q-Hamiltonian setting. 
\begin{lemma}\label{lem:nw}
The distributions $\n_0$ and $\ti{\n}_0$ are related as follows,  
$$ \n_0(\mu)= e^{\tpi\l\mu,q'(\f{1}{\tpi}\f{\p}{\p\nu})\r}
\det(p''( \f{1}{\tpi}\f{\p}{\p\nu}))\,\ti{\n}_0(\nu)\Big|_{\nu=\mu}.
$$
In particular, 
$$ \n_0(0)=\det(p'')\big(\f{1}{\tpi}\f{\p}{\p\nu}\big)\,\ti{\n}_0(\nu)\Big|_{\nu=0}
$$
\end{lemma} 
\begin{proof}
It suffices to prove this for the case that $\phi_0$ is compactly supported. 
The calculation 
$$ \l\n_0,e^{-\tpi\l\cdot,p'(\xi)\r}\r
=\l\phi_0(\xi),e^{-\tpi\l\cdot,p'(\xi)\r}\r
=\l\ti{\phi}_0(\xi),e^{-\tpi\l\cdot,\xi\r}\r=
\l\ti{\n}_0,e^{-\tpi\l\cdot,\xi\r}\r.$$
shows that the Fourier coefficients of $\n_0$ and $\ti{\n}_0$ are
related by a change of variables.  Hence, the Lemma boils down to a
description of the Fourier transform of the change-of-variables operator 
$\mf{A}_\g$, which is accomplished in Appendix \ref{app:fcov}.
\end{proof}

\begin{proof}[Proof of Theorem \ref{th:interpret1}]
With our preparations, the proof has now become a fairly straightforward 
extension of the proof of Theorem \ref{th:interpret}. We may assume that 
$V$ consists of $\beta$-regular values of $\Phi$ and is the diffeomorphic 
image under $\exp$ of an invariant open subset $V_0\subset\g$.
Let $(M,\om_0,\Phi_0)$ be as in the proof of Theorem \ref{th:interpret}. 
We use the currents $\phi=\Phi_*(\beta e^{\tpi\om^p})$ and 
$\phi_0=(\Phi_0)_*(\beta e^{\tpi\om^p_0})$ to define 
$\n^\beta\in\D'(V)^G[[\delta]]$ and $\n_0^\beta \in \D'(V_0)^G[[\delta]]$. 
>From Lemma \ref{lem:duflo2}, we see that $\n^\beta(e)=\n_0^\beta(0)$. 
To identify $\n_0^\beta(0)$ 
it is better to work with $\ti{\n}^\beta_0$, defined by the 
$\d_G-\tpi\eta_\g$-closed form
$$\ti{\phi}_0(\xi)=(\Phi_0)_*(\beta e^{\tpi (\om^p_0-\l\Phi_0, q'(\cdot))})
=\Phi_0^* (\beta e^{\tpi (\om^q_0-\l\Phi_0, q'(\cdot))}
e^{\tpi \om_0})
\in \ca{C}_G(V_0)[[\delta]]. 
$$
Indeed, $\ti{\n}^\beta_0$
is just the DH-distribution
for the $\d_G$-closed equivariant form $\beta e^{\tpi (\om^q_0-\l\Phi_0,
q'(\cdot))}$, hence its interpretation in terms of
intersection pairings is given by Proposition
\ref{prop:nb0}\eqref{it:interpret}.  According to Lemma \ref{lem:nw},
the value of $\n^\beta_0$ at $0$ is equal to the value of $(\det
p'')(\f{1}{\tpi}\f{\p}{\p\mu})\ \ti{\n}^\beta_0$ at $0$.  By 
Proposition \ref{prop:nb0}\eqref{it:polynomial},
$$ (\det p'')\big(\f{1}{\tpi}\f{\p}{\p\mu}\big)\ \ti{\n}_0^\beta=
\ti{\n}^{\det(p'')\beta}_0,$$
so we find that
$\n^\beta(e)=\n^{\beta}_0(0)=\ti{\n}^{\det(p'')\beta}_0(0)$ is given
by
$$ 
(\tpi)^{\dim G}\vol_{G}\sum_j \f{1}{k_j} 
\int_{(M_0\qu G)_j} \big(\det(p'')\,\beta 
e^{\tpi (\om^q_{0}-\l\Phi_0,q'(\cdot)\r)}
\big)_{\on{red}}
e^{\tpi (\om_0)_{\on{red}}}.
$$
Finally, observe that under the natural identification 
$M\qu G=M_0\qu G$, 
$$(\om^q_{0}-\l\Phi_0,q'(\cdot)\r)_{\on{red}}=\om^q_{\on{red}},\
(\om_0)_{\on{red}}=\om_{\on{red}}.$$
\end{proof}

\section{Proof of the Witten formulas}\label{sec:witproof}
In this Section we will apply Theorem \ref{th:higher} to our main
examples. In each case, we first describe the Fourier coefficients of the
distribution $\n^\beta\in \D'(G)[[\delta]]$, by working out the fixed
point contributions
%
from Proposition \ref{prop:locform}. Throughout this Section we assume that 
$G$ is simply connected and simple. Thus in particular, $\rho$ is a weight.
It is easy to generalize the discussion to general semi-simple compact groups, 
by passing to covers. 

\subsection{Conjugacy classes}\label{subsec:conclass}
Any conjugacy class $\Co\subset G$ has a unique structure of a
q-Hamiltonian $G$-space, in such a way that the moment map is the
inclusion $\Phi:\,\Co\hra G$. The 2-form $\om$ given by the formula
\cite{al:mom}
$$ \om(\xi_\Co,\xi'_\Co)|_g=\hh \xi\cdot(\Ad_g-\Ad_{g^{-1}})\xi'.$$
More generally, if $p\in \on{Pol}(\g)^G$ is an invariant polynomial,
the pull-back $\Phi^*\eta^p_G$ is exact since the equivariant
cohomology of a homogeneous space $G/K$ lives in even degree. In fact, 
evaluation at the base point $eK\in G/K$ is a homotopy
equivalence, 
\begin{equation}\label{eq:basept}
\Om_G^\bullet(\Co)\to \Om_K^\bullet(\pt)=\on{Pol}^{2\bullet}(\k)^K,
\end{equation}
with homotopy inverse $\Om_K(\pt)\to
\Om_G(\Co)$ given by induction.  Hence, up to equivariant coboundaries
we can fix a primitive $\om^p\in\Om_G(\Co)$ of $\Phi^*\eta^p_G$ by
requiring that it lies in the kernel of \eqref{eq:basept}. 
Now let $p\in \on{Pol}(\g)[[\delta]]$ have the form \eqref{eq:pform},
and take $\om^p\in\Om_G(\Co)[[\delta]]$ as in \eqref{eq:omp},
normalized to lie in the kernel of \eqref{eq:basept}.
\begin{theorem}
Let $\beta\in\Om_G(\Co)$ be the equivariant cocycle defined by induction 
from an invariant polynomial $Q\in \on{Pol}(\k)^K$. Then the 
Fourier coefficients 
of the corresponding distribution $\n^\beta\in\D'(G)[[\delta]]$ are given 
by, 
$$ \l\n^\beta,\ol{\chi_\lambda}\r=(\tpi)^{\dim(\Co)/2}
Q(\xi) \Big(\f{\det p''(\xi)}{\det p_\k''(\xi)}\Big)^{1/2}
\ol{\chi_\lambda(\Co)}\,\Vol(\Co).$$
with $p'(\xi)=\lambda+\rho$. Here $\chi_\lambda(\Co)$ denotes the
value of the character $\chi_\lambda$ on the conjugacy class, and
$\Vol(\Co)$ is given by \eqref{eq:covol}.
\end{theorem}
\begin{proof}
Let $\mu\in \Alc$ be the unique point in alcove with
$\exp(\mu)\in\Co$, as in Example \ref{ex:conclass}. The centralizer
$K=G_{\exp\mu}$ is a connected subgroup containing the maximal torus
$T$. We denote by $W_K\subset W$ the Weyl group of $K$, by
$\mf{R}_{K,+}\subset \mf{R}_+$ the set of positive roots, and by
$\rho_K$ its half-sum. $\rho_K$ need not be a weight of $K$, in
general.  On the other hand, $2\mu\cdot\rho_K$ is an integer: This
follows since $2\rho_K$ is in the weight lattice for the adjoint group
$K/Z(K)$, while $\mu$ is in the integral lattice for this group (since
$\exp_K(\mu)\in Z(K)$). The orientation on $\Co$ differs from the
orientation on the homogeneous space $G/K$ by a sign,
$(-1)^{2\mu\cdot\rho_K}.$

The set of fixed points for the vector field $(\lambda+\rho)_{\Co}$ is
the Weyl group orbit $W.\exp\mu= \Co\cap T$. It is thus parametrized
by $W/W_K$. Consider the fixed point $F=\{\exp\mu\}$.  The space
$\nu_F=T_{\exp\mu}\Co$ is isomorphic to $\k^\perp\subset \g$ as a
$K$-module; hence the equivariant Euler form is
$$ \Eul(\nu_F,\xi)=(-1)^{2\mu\cdot\rho_K} \prod_{\alpha\in\mf{R}_+\backslash \mf{R}_{K,+}} (-\alpha)\cdot\xi,\ \ \xi\in\t$$
By equivariance, $\Eul(\nu_{w.F},\xi)=\Eul(\nu_F,w\xi)$. 
We therefore obtain, 
\beq \f{\l\n^\beta,\ol{\chi_\lambda}\r}{\dim V_\lambda}
&=&(-1)^{2\mu\cdot\rho_K} Q(\xi) \sum_{w\in W/W_K}
\f{\ol{\eps_{w.(\lambda+\rho)}}(g)}{\prod_{\alpha\in\mf{R}_+\backslash
\mf{R}_{K,+}} (-\alpha)\cdot w\xi}\\ &=&(-1)^{2\mu\cdot\rho_K} Q(\xi)
\Big(\f{\det p''(\xi)}{\det p_\k''(\xi)}\Big)^{1/2} \sum_{w\in W/W_K}
\f{\ol{\eps_{w.(\lambda+\rho)}}(g)}{\prod_{\alpha\in\mf{R}_+\backslash
\mf{R}_{K,+}} (-\alpha)\cdot w(\lambda+\rho)}\eeq
where we have used Lemma \ref{lem:alemma} to write
$$ 
\prod_{\alpha\in\mf{R}_+\backslash \mf{R}_{K,+}} 
\f{\alpha\cdot p'(\xi)}{\alpha\cdot \xi}
=\Big(\f{\det p''(\xi)}{\det p_\k''(\xi)}\Big)^{1/2}
$$
The proof is complete by the following Lemma. 
\end{proof}

\begin{lemma}
Let $K$ be the centralizer of $g=\exp\mu$, with $\mu\in\Alc$. 
Let $\Co=G.\exp\mu$ be the corresponding conjugacy class. 
The following formula holds for all $\lambda\in \Lambda^*_+$: 
$$ \f{\chi_\lambda(\Co)}{\dim V_\lambda}=
\f{(-1)^{2\mu\cdot\rho_K}}{ \Vol(\Co)}\,
(\tpi)^{-\dim(\Co)/2}
\sum_{w\in
W/W_K}\f{ \eps_{w.(\lambda+\rho)}(g)}{\prod_{\alpha\in\mf{R}_+\backslash
\mf{R}_{K,+}} \alpha\cdot w(\lambda+\rho)},$$
with $\Vol(\Co)$ as given in \eqref{eq:covol}.
\end{lemma}
\begin{proof}
Let $\ti{K}$ denotes any connected finite cover of $K$ for which $\rho_K$ is in
the weight lattice, and $\chi_\nu^{\ti{K}}\in C^\infty(K)$ denotes the
irreducible character corresponding to a weight $\nu\in\Lambda^*_{\ti{K},+}$.
By two applications of the Weyl character formula, 
\beq \chi_\lambda(\Co)&=&\sum_{w\in W/W_K} (-1)^{l(w)}
\f{\chi^{\ti{K}}_{w(\lambda+\rho)-\rho_K}(\exp_{\ti{K}}\mu)}{\prod_{\alpha\in\mf{R}_+\backslash \mf{R}_{K,+}}2i\sin(\pi\alpha\cdot\mu) }
\\&=&\f{\vol_{G/K}}{\Vol(\Co)} i^{-\dim(\Co)/2} 
\sum_{w\in W/W_K} (-1)^{l(w)}
\chi^{\ti{K}}_{w(\lambda+\rho)-\rho_K}(\exp_{\ti{K}}\mu).
\eeq
Since $\exp_{\ti{K}}\mu$ is in the center of $\ti{K}$, 
$$\chi_\nu^{\ti{K}}(\exp_{\ti{K}}\mu)= \dim V_\nu^{\ti{K}}\ e^{\tpi
\nu\cdot\mu}$$ 
for all $\nu\in \Lambda^*_{\ti{K},+}$. In our case
$\nu=w(\lambda+\rho)-\rho_K$, and therefore $e^{\tpi \nu\cdot\mu}=
(-1)^{2\rho_K\cdot\mu} 
e^{\tpi w(\lambda+\rho)\cdot \mu}$, while the dimension of
$V_\nu^{\ti{K}}$ can be re-written by means of the Weyl dimension formula:
\beq \dim V_{w(\lambda+\rho)-\rho_K}^{\ti{K}}&=&
\vol_{K/T}\prod_{\alpha\in \mf{R}_{K,+}} 2\pi \alpha\cdot
w(\lambda+\rho)\\ &=&\f{(-1)^{l(w)}\dim V_\lambda
}{\vol_{G/K}\prod_{\alpha\in \mf{R}_+\backslash \mf{R}_{K,+}} 2\pi
\alpha\cdot w(\lambda+\rho)}. \eeq
Putting all this together, we obtain the expression for 
$\chi_\lambda(\Co)$ as given in the Lemma. 
\end{proof}

\subsection{Fusion}
Suppose $(M,\om,\Phi)$ is a q-Hamiltonian $G^2$-space. 
That is, $\Phi=(\Phi_0,\Phi_1):\,M\to G^2$ is a $G^2$-equivariant map, 
and the 2-form $\om$ satisfies the condition 
$\d_{G^2}\om=\Phi_0^*\eta_G+\Phi_1^*\eta_G$. Then, as explained 
in \cite{al:mom} the space $M$ with diagonal $G$-action becomes 
a q-Hamiltonian $G$-space $(M,\om_{\on{fus}},\Phi_{\on{fus}})$ 
with moment map $\Phi_{\on{fus}}=\Phi_0\Phi_1$ (pointwise product) and
2-form 
$$\om_{\on{fus}}=\om +\Phi^*\varphi,$$
where $\varphi=\hh \pr_1^*\theta^L\cdot \pr_2^*\theta^R\in 
\Om^2(G\times G)$. 
In the appendix, we explain how the Bott-Shulman machinery defines higher 
analogues $\varphi^p\in \Om_G^{2\bullet-2}(G\times G)$ of this form, for any 
$p\in \on{Pol}^\bullet(\g)^G$, with $\varphi^p=\varphi$ for 
$p(\xi)=\hh||\xi||^2$. Some basic properties of these forms are: 
\begin{enumerate}

\item\label{item:a} 
The form degree $0$ part $(\varphi^p)^{[0]}$ vanishes. 
\item\label{item:d}
Letting $\on{Mult}:\,G\times G\to G$ denote the group multiplication, 
$$ \on{Mult}^*\eta^p_G=\pr_1^*\eta^p_G+\pr_2^*\eta^p_G+\d_G\varphi^p.$$
\item\label{item:b} 
The equivariant forms on $G\times G\times G$,
$$(\pr_1\times \pr_2)^*\varphi^p+
(\on{Mult}\circ (\pr_1\times \pr_2)\times \pr_3)^*\varphi^p$$
and 
$$(\pr_2\times\pr_3)^*\varphi^p+(\pr_1\times (\on{Mult}\circ (\pr_2\times \pr_3))^*\varphi^p$$
differ by a $\d_G$-coboundary. For $p(\xi)=\hh||\xi||^2$ they are in 
fact equal. 
\end{enumerate}

Hence, given an equivariant form $\om^p\in \Om_{G^2}(M)$ 
with $\d_{G^2}\om^p=\Phi_0^*\eta_G^p+\Phi_1^*\eta_G^p$, the form 
$$\om^p_{\on{fus}}=\om^p+\Phi^*\varphi^p$$
has the property $\d_G\om^p_{\on{fus}}=\Phi_{\on{fus}}^*\eta^p_G$.

As an application, suppose $\Co_1,\ldots,\Co_r$ are conjugacy classes,
with moment maps $\Phi_j$ the inclusion and with given higher
q-Hamiltonian forms $\om^p_j$, normalized as in \ref{subsec:conclass}. 
Then the product $\Co_1\times
\cdots\times \Co_r$ becomes a q-Hamiltonian $G$-space, with moment map
the product $\Phi_1\cdots\Phi_r$, and with a higher q-Hamiltonian form
$\om^p$ obtained by iterated fusion. Note that $\om^p$
is not canonically defined since the fusion procedure is not
associative (except for $p(\xi)=\hh||\xi||^2$). However,
\eqref{item:b} above shows that any two forms obtained in this way
differ by a $\d_G$-coboundary.

The fixed point set for the action of $(\lambda+\rho)_M$ on
$M=\Co_1\times \cdots\times \Co_r$ are simply products of the
fixed point sets for each $\Co_j$.  Since they are 0-dimensional, the
terms $\Phi^*\varphi^p$ vanish if pulled back to the fixed point sets,
using \eqref{item:a}. Hence, the fixed point contribution for the product 
is just the product of fixed point contributions from each factor.

\subsection{The double}
The {\em double} is the q-Hamiltonian $G^2$-space
$(D(G),\om,(\Phi_0,\Phi_1))$ where $D(G)=G\times G$, with group action
$$ (g_0,g_1).(h,k)=(g_0 h g_1^{-1},\,g_1 k g_1^{-1})$$ 
and moment map components $\Phi_0(h,k)=k$, 
$\Phi_1(h,k)=\Ad_h(k^{-1})$. The 2-form $\om$ satisfying 
$\d_{G^2}\om=\Phi_0^*\eta_G+\Phi_1^*\eta_G$ arises naturally 
in the Bott-Shulman construction (cf. Appendix \ref{subsec:G}), and 
equals the form denote $\lambda$ in \eqref{eq:lambda}.  
The explicit form of $\om$ is not important for what follows, 
except for the fact that its pull-back to $D(T)=T\times T\subset G\times G$ 
is the standard symplectic form on $D(T)=T\times T$: 
$$ \iota_{D(T)}^*\om= h^*\theta_T\cdot k^*\theta_T.$$
More generally, for any invariant polynomial $p\in \on{Pol}(\g)^G$,
the Bott-Shulman construction defines forms $\om^p\in \Om_{G^2}(D(G))$
(see Appendix \ref{subsec:G}, where these forms are denoted
$\lambda^p$), with the property
$\d_{G^2}\om^p=\Phi_0^*\eta_G^p+\Phi_1^*\eta_G^p$. The pull-back
of these forms to $D(T)$ is given by
$$ \iota_{T\times T}^*\om^p(\xi_0,\xi_1)=\sum_{jk} A_{rs}(\xi_0,\xi_1)
h^*\theta_T^r
k^*\theta_T^s$$
with $A_{rs}(\xi_0,\xi_1)=\int_0^1 \d t\, p_\t''(\xi)_{rs}((1-t)\xi_0+t\xi_1)$. 
Here $p_\t$ is the restriction of $p$ to $\t$, $p_\t''(\xi)_{jk}$ is the 
matrix of second derivatives in a given orthonormal basis of $\t$, and 
$\theta_T^j$ are the corresponding components of the Maurer-Cartan form. 

The fused double $\ti{D}(G)=D(G)_{\on{fus}}$ is a q-Hamiltonian $G$-space 
with moment map $\Phi_{\on{fus}}=\Phi_0\Phi_1$ given by the group commutator
$(h,k)\mapsto [k,h]=khk^{-1}h^{-1}$. As explained above, the forms $\om^p$ 
gives rise to higher q-Hamiltonian forms $\om^p_{\on{fus}}$ on 
$\ti{D}(G)$. On $\ti{D}(T)\subset \ti{D}(G)$, the 'fusion term'
$\Phi^*\varphi^p$ vanishes, as one can see from Lemma 
\ref{lem:varphiab}. Hence, the pull-back of $\om^p_{\on{fus}}$ to $\ti{D}(T)$ 
is simply 
$$ \iota_{T\times T}^*\om^p_{\on{fus}}(\xi)=\sum_{jk} 
p_\t''(\xi)_{rs} h^*\theta_T^r
k^*\theta_T^s$$

Dropping the subscript '$\on{fus}$', consider now the fixed point contributions 
for $(\ti{D}(G),\Phi,\om)$, for $p$ be as in \eqref{eq:pform}, and 
$\om^p$ obtained by fusion as explained above. 
Since $G$ acts by conjugation on each factor in $\ti{D}(G)=G\times G$,  
the fixed point set for $(\lambda+\rho)_M$ is just $F=\ti{D}(T)$. 
On $F$ the moment map becomes trivial, and the normal bundle is 
$T$-equivariantly isomorphic to ${\mf{p}}\oplus {\mf{p}}$. 
The equivariant Euler form is (cohomologous to)
$$\Eul(\nu_F,\xi)=(-1)^{n_+}\,(\prod_{\alpha\in\mf{R}_+}\alpha\cdot \xi)^2$$
In the simplest case $\beta=1$, the fixed point contributions read, 
\beq  \f{\l\n,\ol{\chi_\lambda}\r}{\dim V_\lambda}
&=&(-1)^{n_+} \f{\int_F e^{\tpi \om^p(\xi)}}{(\prod_{\alpha\in\mf{R}_+}\alpha\cdot \xi)^2}
\\
&=& (-1)^{n_+} (\tpi)^{\dim T}\, (\vol_T)^2\,
\f{\det p_\t''(\xi)}{(\prod_{\alpha\in\mf{R}_+}\alpha\cdot \xi)^2}\\
&=& (-1)^{n_+}  (\tpi)^{\dim T}\, (\vol_T)^2\,
\f{\det p''(\xi)}{(\prod_{\alpha\in\mf{R}_+}\alpha\cdot (\lambda+\rho))^2}\\
&=& (-1)^{n_+}  (\tpi)^{\dim T}\, \f{(\vol_T)^2}{(\dim V_\lambda)^2}\,
\f{\det p''(\xi)}{(\prod_{\alpha\in\mf{R}_+}\alpha\cdot \rho)^2}\\
&=& (\tpi)^{\dim G} 
\f{(\vol_G)^2 \det p''(\xi)}{(\dim V_\lambda)^2}.
\eeq
We now consider more general $\beta$. The equivariant cohomology algebra of 
$G\times G$ is a tensor product (as modules over $H_G(\pt)$) of the 
cohomology algebras of the two factors
$$ H_G(G\times G)=H_G(G)\otimes_{H_G(\pt)}H_G(G).$$
Let $p_1,\ldots,p_N$ be a set of generators for the ring of invariant
polynomials, and let $\eta^i_G=\eta^{p_i}_G$ be the corresponding
generators of $H_G(G)$. Following \cite{wi:tw} and \cite[Section
6]{al:gr}, set
$$ \beta(\xi)=\exp(\sum_{j=1}^N \sig_j p_j(\xi)+\sum_{i=1}^N \eps_i^{(1)} h^*\eta_G^{p_i}(\xi)+\eps_i^{(2)} 
k^*\eta^{p_i}_G(\xi)).$$
where $\sig_i$ are formal even parameters, and $\eps_i^{(1)},
\eps_i^{(2)}$ are {\em odd} parameters (anti-commuting with each other
but also with odd differential forms). To simplify notation, denote
$$P^{(1)}(\xi)=\sum_{i=1}^N \eps_i^{(1)} p_i(\xi),\ \ \
P^{(2)}(\xi)=\sum_{i=1}^N \eps_i^{(2)} p_i(\xi),\ \ \ R(\xi)=
\sum_{i=1}^N \sig_i p_i(\xi)$$
Thus
$$ \beta=\exp(R+h^*\eta_G^{P^{(1)}}+k^*\eta_G^{P^{(2)}}).$$
The integral over $F=T^2$ now becomes
$$ \int_{T^2} \beta(\xi) e^{\tpi \om^p(\xi)}
=e^{R(\xi)} \int_{T^2} e^{\tpi (p_\t''(\xi)h^*\theta^T)\cdot k^*\theta^T)
+(P^{(1)})'(\xi)\cdot h^*\theta^T + (P^{(2)})'(\xi)\cdot k^*\theta^T}.$$
This integral is solved by completion of the square, 
writing the exponent as 
$$ \big(\tpi\ p_\t''(\xi)h^*\theta^T+Q'(\xi)\big)\cdot 
\big(k^*\theta^T-\f{1}{\tpi} p_\t''(\xi)^{-1}P'(\xi)\big)\,
-\, \f{1}{\tpi} p_\t''(\xi)^{-1}(P^{(1)})'(\xi) \cdot (P^{(2)})'(\xi).$$
This yields,  
$$ \l\n^\beta,\ol{\chi_\lambda}\r=(\tpi)^{\dim G} 
\f{(\vol_G)^2 \det p''(\xi)}{(\dim V_\lambda)^2}
e^{R(\xi)-\, \f{1}{\tpi} p_\t''(\xi)^{-1}(P^{(1)})'(\xi) \cdot 
(P^{(2)})'(\xi)}
$$
with $p'(\xi)=\lambda+\rho$. 

\subsection{Moduli spaces of flat bundles}
By fusion of $s$ copies of the double $G^2$ and $r$ conjugacy classes
$\Co_1,\ldots,\Co_r$, the space $M=G^{2s}\times \Co_1\times\cdots
\times\Co_r$ considered in \eqref{eq:spaceM} acquires the structure of
a q-Hamiltonian $G$-space, with moment map \eqref{eq:mommap}.  For any
$p\in \on{Pol}(\g)^G$, this space carries a higher q-Hamiltonian form
$\om^p$, which is canonically defined up to a $\d_G$-coboundary.

The equivariant cohomology ring of $M$ is simply the tensor product 
(over $H_G(\pt)$) of $2s$ factors of $H_G(G)$ with the cohomology 
rings of the conjugacy classes $\Co_l=G/K_l$, i.e. 
$\on{Pol}(\k_l)^{K_l}$. Let $R=\sum \sig_i p_i$ be as above, and write
$$
P^{(1)}_j=\sum_{i=1}^N \eps_{ij}^{(1)} p_i,\ \ \
P^{(2)}_j=\sum_{i=1}^N \eps_{ij}^{(2)} p_i$$
where $\eps_{ij}^{(1)}, \eps_{ij}^{(2)}$ are odd parameters, and the
index $j$ stands for the $G^2$-factor.
Let $Q_l\in \on{Pol}(\k_l)^{K_l}$ be 
invariant polynomials, with corresponding cocycles 
$\beta_l\in\Om_G(\Co_l)$, and let 
$$ \beta=\exp\Big(R+\sum_{j=1}^s (h_j^*\eta_G^{P^{(1)}_j}+
k_j^*\eta_G^{P^{(2)}_j})\Big)\prod_{l=1}^r \beta_l$$
Consider now the fixed point contributions.  The fixed point sets $F$
for $\lambda+\rho$ is just the product of the fixed point sets for the
factors $G^2$ resp. $\Co_l$. In particular, each of the $r+s$ moment
maps is {\em constant} on $F$. This implies that the `fusion terms',
i.e.  terms involving $\varphi^p$, all vanish if pulled back to
$F$. Hence, the fixed point contribution from any such $F$ is simply
the product of the fixed point contributions from the factors.
We therefore obtain the following result for the 
Fourier coefficients, 
$$ \l\n^\beta,\ol{\chi_\lambda}\r=(\tpi)^{\f{\dim M}{2}} (\vol_G)^{2s}
\Big(\f{{\det}^{1/2} p''(\xi)}{\dim V_\lambda}\Big)^{2s+r}\,e^{\ti{R}(\xi)}\,\prod_{l=1}^r 
 \f{Q_l(\xi)\Vol(\Co_l)\ol{\chi_\lambda(\Co_l)}}{{\det}^{1/2} p_{\k_l}''(\xi)  }
$$
where 
$$\ti{R}(\xi)=R(\xi)-\, \f{1}{\tpi} \sum_j
p_\t''(\xi)^{-1}(P_j^{(1)})'(\xi) \cdot (P_j^{(2)})'(\xi)
$$
We finally state the resulting Witten formula for intersection
pairings.  Since we assume $G$ is simple and simply connected, the
level set $\Phinv(e)\subset M$ is connected, and so is the moduli space 
$\M=M\qu G$. If
$e$ is a regular value, the action of $G$ on $\Phinv(e)$ is locally
free. We assume that the generic stabilizer for the $G$-action on
$\Phinv(e)$ is equal to the center $Z(G)$. This is automatic for $s\ge
2$, or for $2s+r\ge 3$ and sufficiently 'generic' conjugacy classes.
(See \cite{al:du} for discussion.) Write $\mf{f}^p=\om^p_{\on{red}}$ 
since these generalize the Atiyah-Bott classes of 'type $\mf{f}$', 
and let $\beta$ be given as above. Then

\beq \lefteqn{\f{1}{(\tpi)^{\dim \M/2}} \int_{\M} \beta_{\on{red}}
e^{\tpi \mf{f}^p}}\\ &=&\#Z(G) (\vol_G)^{2s-2}
\sum_{\lambda\in\Lambda^*_+} e^{\ti{R}(\xi)} \Big( \f{{\det}^{1/2}
p''(\xi)}{\dim V_\lambda}\Big)^{2s+r-2}
\prod_{l=1}^r\Big(\f{Q_l(\xi)\Vol(\Co_l)\ol{\chi_\lambda(\Co_l)}}{{\det}^{1/2}
p_{\k_l}''(\xi)}\Big)\eeq
\begin{appendix}

\section{Formal change of variables}\label{app:fcov}
A {\em formal diffeomorphism} of a manifold $X$, 
is an invertible elements of the algebra 
$C^\infty(X,X)[[\delta]]$ where $\delta_j$ are given formal parameters. 
If $X=V$ is a vector space, we can consider the smaller group 
of {\em polynomial} formal diffeomorphisms, consisting of invertible 
elements $P=\sum_I \delta^I P_I$ (using multi-index notation) 
of the algebra $\on{Pol}(V,V)[[\delta]]$, where 
$\on{Pol}(V,V)$ are $V$-valued polynomials on $V$. 

Suppose $V$ carries a scalar product $\cdot$, and let
$\F_V^{-1}:\,\D'(V)_{\on{comp}}\to C^\infty(V)$ denote the inverse
Fourier transform, defined by
$$ (\F_V^{-1} \n)(\xi)=\l\n,e^{-\tpi\l\cdot,\xi\r}\r,\ \ \xi\in V.$$
\begin{proposition}\label{prop:formal}
Suppose $P=\sum_I \delta^I P_I$ is a formal diffeomorphism of $V$, 
such that all $P_I$ are polynomials and $P_\emptyset=\on{Id}_V$. Let 
$\n_1,\n_2\in \D'(V)_{\on{comp}}[[\delta]]$, with 
$$ (\F_V^{-1}  \n_1)(\xi)=(\F_V^{-1}  \n_2)(P^{-1}(\xi))$$
as elements of $C^\infty(V)[[\delta]]$. Then 
$$ \n_1(\mu)=e^{\tpi \mu\cdot Q(\f{1}{\tpi}\f{\p}{\p \nu})}\,
\det\big(P'(\f{1}{\tpi}\f{\p}{\p \nu})\big)\,
\n_2(\nu)\Big|_{\nu=\mu}
$$
where $Q(\xi)=P(\xi)-\xi$. 
\end{proposition}
\begin{proof}
Results of this type are well-known from the theory of Fourier integral 
operators -- see e.g. \cite{du:fo}. It is enough to consider the case 
$\n_i$ smooth, so that  $f_i(\xi)=(\F_V^{-1} \n_i)(\xi)$ are 
rapidly decreasing functions of $\xi$. 
The desired identity follows from the calculation, 
\beq \n_1(\mu)
&=& \int_V f_2(P^{-1}(\xi)) e^{\tpi \mu\cdot\xi}\d \xi\\ &=& \int_V
f_2(\zeta)\,e^{\tpi \mu\cdot P(\zeta)}(\det P'(\zeta))\d \zeta\\ &=&
\int_V e^{\tpi \mu\cdot Q(\zeta)}(\det P'(\zeta))\, f_2(\zeta)\,e^{\tpi
\mu\cdot\zeta}\d \zeta.\eeq
(To justify the second equality, use Borel summation to temporarily replace
$P$ by a genuine function $V\to V$, depending on $\delta_j$ as parameters.) 
\end{proof}

\section{Equivariant de Rham theory}
\subsection{The Cartan model of equivariant cohomology}
\label{app:carmod}
Let $G$ be a compact Lie group acting smoothly on a manifold $M$. That
is, we are given a group homomorphism $G\to \on{Diff}(M),\,g\mapsto
\ca{A}_g$, such that the action map $G\times M\to M,\ (g,x)\mapsto
g.x=\ca{A}_g(x)$ is smooth. For any Lie algebra element $\xi\in \g$,
the corresponding generating vector field $\xi_M$ is the derivation of
$C^\infty(M)$ given by $\xi_M(f)=\f{d}{d
t}|_{t=0}\ca{A}_{\exp(-t\xi)}^*f$.  For each $\xi\in\g$, define an odd
derivation of $\Om(M)$ by
$$ \d_\xi=\d-\iota(\xi_M).$$
The Cartan complex of equivariant differential forms is the graded
algebra $\Om_G^\bullet(M)=(\on{Pol}(\g)\otimes\Om(M))^G$ of
$G$-equivariant polynomial maps $\beta:\,\g\to \Om(M)$, equipped with
the equivariant differential
$$ (\d_G\beta)(\xi)=\d_\xi\beta(\xi),$$
and with grading given by the differential form degree plus twice the
polynomial degree. Its cohomology algebra coincides with the
equivariant cohomology $H_G^\bullet(M)=H^\bullet(EG\times_G M)$.

\subsection{The Cartan map}\label{app:car}
Let $G$ be a compact Lie group and $\pi:\,P\to B$ a principal $G$-bundle. 
Then the pull-back map $\pi^*:\,\Om(B)\to \Om(P)_{\on{basic}}\subset \Om_G(P)$ 
is a quasi-isomorphism. If $\theta\in \Om^1(P)\otimes\g$ is a principal
connection, one has an explicit homotopy inverse
$$\on{Car}^\theta:\ \Om_G(P)\to \Om(P)_{G-\on{basic}}\cong \Om(B),$$ 
known as the {\em Cartan map}. The definition of this map is as follows: Let
$\on{Hor}^\theta:\,\Om^\bullet(P)\to \Om^\bullet(P)$ denote the
horizontal projection, and let
$F^\theta=\d\theta+\hh [\theta,\theta]\in \Om^2(P)\otimes\g$ be the
curvature. There is a unique algebra homomorphism 
$$ \on{Pol}^\bullet(\g)\to \Om^{2\bullet}(P),\ p\mapsto p(F^\theta)$$
given on linear polynomials $\mu\in \g^*$ by $\mu\mapsto \l\mu,F^\theta\r$.
Tensoring with $\Om(P)$, this yields an algebra homomorphism 
$$ \Om_G(P)=(\on{Pol}(\g)\otimes \Om(P))^G\to \Om(P)^G,\ \beta\mapsto
\beta(F^\theta).$$
The Cartan map is defined by 
$$ \on{Car}^\theta(\beta)=\on{Hor}^\theta(\beta(F^\theta)).$$
It was proved by Cartan that $\on{Car}^\theta$ is a chain map,
inducing the inverse map $(\pi^*)^{-1}$ in cohomology. For a nice
proof of Cartan's theorem, showing in particular that
$\on{Car}^\theta$ and $\pi^*$ are homotopy inverses, see Nicolaescu
\cite{ni:ca}. The proof carries over to the case that the $G$-action
on $P$ is not free but only {\em locally} free, i.e. has finite
stabilizers, and it also generalizes to the case that $P\to B$ 
is an $L$-equivariant principal $G$-bundle with a $L$-invariant 
connection, where $L$ is a second Lie group. 

One application of the Cartan map is {\em induction}. Suppose $G$ is a
compact connected Lie group, and $K$ a maximal rank subgroup. Let the
principal $K$-bundle $G\to G/K$ be equipped with the unique
$G$-invariant connection.  Given a $K$-manifold $Y$, let the principal
$K$-bundle $G\times Y\to G\times_K Y$ carry the pull-back connection.
The {\em induction map}
$$ \Ind_K^G:\,\Om_K(Y)\to \Om_G(G\times_K Y),$$
is a homotopy equivalence, given as the pull-back map 
$\Om_K(Y)\to \Om_{G\times K}(G\times Y)$ followed by 
the Cartan map. A homotopy inverse 
$\Om_G(G\times_K Y)\to \Om_K(Y)$ is given as pull-back to 
$Y\subset G\times_K Y$. 

\begin{example}\label{ex:hom}
In the special case $Y=\pt$, we obtain homotopy inverses
$$ \Om_K(\pt)\to \Om_G(G/K),\ \ \Om_G(G/K)\to \Om_K(\pt)$$
where the first map is the induction map and the second map is
pull-back to the base point $eK\in G/K$.
\end{example}

\subsection{Equivariant Thom form, equivariant Euler form}
Let $\pi:\,E\to B$ be a $G$-equivariant real vector bundle of even
rank, with a fiberwise orientation. An equivariant Thom form for $E$
is an equivariant form $\tau_E\in\Om_G(E)$, compactly supported in
fiber directions, with fiber integral $\pi_*\tau_E=1$.  Any two Thom
forms differ by the equivariant coboundary of a form with fiberwise
compact support. Assume the base $B$ is oriented, and give $E$ the
product orientation. Then
\begin{equation}\label{eq:thomproperty} 
\int_E \tau_E(\xi)\beta =\int_B \beta, 
\end{equation}
for any differential form $\beta\in \Om(E)$ such that $\tau(\xi)\beta$
has compact support and $\d_\xi\beta=0$. The proof boils down to the
fact that if one is working with equivariant currents, the Thom class
is represented by $\iota_*(1)\in \ca{C}_G(E)$, where $\iota:\,B\to E$
is the inclusion of the base.

The pull-back $\on{Eul}(E)=\iota^*\tau_E\in\Om_G(B)$ is
called the equivariant Euler form. Given an invariant Riemannian
metric and compatible connection on $E$, the Mathai-Quillen
construction \cite{ma:th} gives explicit representatives for $\tau_E$,
and therefore of $\on{Eul}(E)$.

An important special case is: $G=T$, $B=\pt$, and $E=V$ a complex vector 
space. Let $a_1,\ldots,a_n\in\Lambda^*$ be the (real) 
weights for the action on $V$. Then 
the equivariant Euler form is simply a polynomial on $\t$: 
$\Eul(V,\xi)=(-1)^n\prod_{j=1}^n \l a_j,\xi\r$. 

\subsection{The Berline-Vergne localization formula}\label{app:bv}
Let $G$ be a compact Lie group and $M$ an oriented $G$-manifold. For
any $\xi\in\g$, consider the derivation $\d_\xi=\d-\iota(\xi_M)$.  Let
$M^\xi$ be the set of zeroes of $\xi_M$, or equivalently the fixed
point set of the 1-parameter subgroup generated by $\xi$.

\begin{theorem}[Berline-Vergne] \cite{be:cl}\label{th:bv} 
Suppose $\alpha\in \Om(M)_{\on{comp}}$ is a compactly supported
differential form with $\d_\xi\alpha=0$.  Let $S\subset M$ be an
embedded $G_\xi$-invariant oriented submanifold of even codimension,
containing the fixed point set $M^\xi$. (We allow for $S$ to consist
of several components of varying dimension.)
Let $\Eul(\nu_S,\cdot)\in \Om_{G_\xi}(S)$ be the $G_\xi$-equivariant
Euler form of the normal bundle of $S$, for some choice of invariant
Euclidean metric and connection.
Then
$$ \int_M \alpha =\int_S \f{\iota_S^*\alpha}{\Eul(\nu_S,\xi)}$$
\end{theorem}

\begin{remarks}
\begin{enumerate}
\item
Suppose $\beta=\alpha(\xi)$ where $\alpha\in\Om_G(M)$ is a
$G$-equivariant cocycle. Then $\d_\xi\beta=0$, and the localization
formula is a version of the localization formula in equivariant
cohomology (see Atiyah-Bott \cite{at:mom}). However, not all
$\d_\xi$-cocycles arise in this way.
\item 
The theorem is usually stated for the special case $S=M^\xi$.  In
fact, the general case may be deduced from this special case, by further
localizing the integral over $S$ to $S'=M^\xi\subset S$. If $\zeta\in \g$
is sufficiently close to $\xi$ and commutes with $\xi$, then also
$S=M^\zeta$ satisfies the conditions of the theorem.
\end{enumerate}
\end{remarks}

\subsection{Equivariant homotopy operators}\label{app:hom}
We will need the following facts about homotopy operators.  For any 
$G$-equivariant vector bundle $\pi:\,E\to B$, let $\on{h}:\,\Om^\bullet(E)\to
\Om^{\bullet-1}(E)$ denote the standard homotopy operator. That is, up
to a sign $\on{h}$ is defined as pull-back under the map $I\times E\to E$
given by scalar multiplication on the fibers, followed by the
push-forward map $(\pr_2)_*:\,\Om^\bullet(I\to E)\to
\Om^{\bullet-1}(E)$.  Since the projection $\pi$ is $G$-equivariant, 
the homotopy operator defines a degree $-1$ operator on $\Om_G(E)$, 
denoted by the same symbol. 

Letting $\iota:\,B\to E$ be the inclusion of the
zero section, the homotopy operator satisfies $\iota^*\circ \on{h}=0$,
$\on{h}\circ \pi^*=0$ and
$$ \d_G \on{h} +\on{h} \d_G=\on{id}-\pi^*\circ \iota^*.$$ 
Another simple fact regarding $\on{h}$ is that if a form on $E$ is
zero along $B\subset E$, then so is its image under $\on{h}$.

\subsection{Equivariant simplicial differential forms}\label{app:bs}
Recall the definition of a simplicial manifold
\cite{du:si,se:cl,mo:no}. For each positive integer $n$ let $[n]$
denote the ordered sequence $\{0,\ldots,n\}$.  A map $f:\,[m]\to[n]$
is called {\em increasing} if $f(i)\ge f(j)$ for $i>j$.  Of particular
interest are the {\em face maps} $\partial^i:\,[n-1]\to [n]$ for
$i=0,\ldots,n$, defined as the unique strictly increasing map 
whose image does not contain $i$.

A simplicial manifold is a contravariant functor from the category of
ordered sequences (with increasing maps as morphisms) into the
category of manifolds. That is, a simplicial manifold $X_\bullet$ is a
sequence of manifolds $(X_n)_{n=0}^\infty$, together with a map
$X(f):\,X_n\to X_m$ for each increasing map $f:\,[m]\to [n]$, such
that $X(\on{id})=\on{id}$ and $X(f\circ g)=X(g)\circ X(f)$. The maps
$\partial_i=X(\partial^i):\,X_n\to X_{n-1}$ are again referred to as
{\em face maps}. A (smooth) simplicial map between simplicial manifolds
$F_\bullet:\, X_\bullet\to X_\bullet'$ is a collection of smooth maps
$F_n:\,X_n\to X'_n$ intertwining the maps $X(f),\,X'(f)$. Any manifold 
$M$ can be viewed as a simplicial manifold $M_\bullet$, 
where all $M_n=M$ and all $X(f)$ are the identity map. Another example 
is the simplicial manifold $E_\bullet M$, where $E_nM=M^{n+1}$ and
$X(f)(x_0,\ldots,x_n)=(x_{f(0)},\ldots,x_{f(n)})$.  

Let $\Del^n\subset \R^{n+1}$ denote the standard $n$-simplex, defined
as the intersection of the positive orthant with the affine hyperplane
$\sum_{i=0}^n t_i=1$. Any increasing map $f:\,[m]\to [n]$ defines a
linear map $\R^{m+1}\to \R^{n+1}$ sending the basis vector $e_i$ to
$e_{f(i)}$. It induces a map $\Delta(f):\,\Del^m\to \Del^n$. 
The geometric realization \cite{mi:ge,se:cl} of a simplicial manifold
$|X|$ is the quotient $|X|=\coprod_n (\Delta^n\times X_n)/\sim$ where
one divides by the equivalence relation generated by
$(\Delta(f)(t),x)\sim (t,X(f)(x))$ for all increasing maps $f$.  

Following Dupont \cite{du:si}, one defines  
a simplicial $r$-form on $X_\bullet$ to be a collection of
$r$-forms $\alpha_n\in \Om^r(\Delta^n\times X_n)$ satisfying relations
\begin{equation}\label{eq:dupont}
(\Delta(f)\times \id)^*\alpha_n =(\id\times X(f))^*\alpha_{m}\ \ 
\end{equation}
for any increasing map $f:\,[m]\to [n]$. Under certain technical
hypothesis (which holds in our examples), the complex
$(\Om^r_{\on{simp}}(X_\bullet),\d)$ computes the cohomology of the
geometric realization with coefficients in $\R$.  Consider on the
other hand the double complex $\Om^{k,l}(X_\bullet):=\Om^l(X_k)$, with
commuting differentials $\d$ and $\delta=\sum_{i=0}^{k+1}
(-1)^i\partial_i^*$, and the corresponding total complex
$\Om^r(X_\bullet) =\bigoplus_{k+l=r}\Om^l(X_k)$ with
differential $\delta+(-1)^k\d$. The maps 
$\Om^r(\Del^n\times X_n)\to \Om^{r-n}(X_n)$ (integration over 
simplices) assemble to a chain map 
$$ \Om^r_{\on{simp}}(X_\bullet)\to \Om^r(X_\bullet),$$
As shown by Dupont \cite{du:si}, this map is a chain 
homotopy equivalence.  

There is a straightforward equivariant extension of these concepts:
Suppose $K_\bullet$ is a simplicial Lie group (i.e. all $K_n$ are Lie
groups and all face maps and degeneracy maps are group homomorphisms).
An action of $K_\bullet$ on $X_\bullet$ is a simplicial map
$K_\bullet\times X_\bullet\to X_\bullet$ given by a $K_n$-action on
$X_n$ in each degree $n$. This means that for any increasing map
$f:\,[m]\to [n]$, the maps $X(f):\,X_n\to X_m$ are equivariant with
respect to the homomorphisms $K(f):\,K_n\to K_m$:
$$ X(f)(k.x)=(K(f)(k)).(X(f)(x)).$$
Thus one obtains pull-back maps in equivariant cohomology, 
$ X(f)^*:\,\Om^r_{K_m}(X_m)\to \Om^r_{K_n}(X_n).
$
We define a space 
$$\Om^r_{K_\bullet}(X_\bullet):= \bigoplus_{n=0}^r \Om^{r-n}_{K_n}(X_n)$$ 
of $K_\bullet$-equivariant forms on $X_\bullet$, with 
equivariant differential $\delta+(-1)^n\d_{K_n}$ on 
$\Om^{r-n}_{K_n}(X_n)$, and define a 
$K_\bullet$-equivariant $r$-form to be a collection of equivariant 
forms $\alpha_n= \Om^r_{K^n}(\Delta^n\times X_n)$
satisfying the compatibility
relations \eqref{eq:dupont}. As before, integration over simplices 
defines a chain equivalence between these two complexes.

Suppose now that $P_\bullet\to X_\bullet$ is a simplicial
$K_\bullet$-equivariant principal $G$-bundle. A $K_\bullet$-invariant
simplicial connection $\sig_\bullet$ is given by a family of
$K_n$-invariant connection forms 
$$\sig_n\in\Om^1(\Del^n\times
P_n)^{K_n}\otimes\g$$ 
satisfying Dupont's compatibility relations.
Given a $G$-manifold $M$, the collection of equivariant Cartan
maps for $\sig_\bullet$, followed by integration over simplices,
defines a chain map\footnote{See \cite{al:cw} for an alternative construction
of such a chain map.}
\begin{equation}\label{eq:car1} 
\Om^\bullet_G(M)\to 
\bigoplus_n \Om^{\bullet-n}_{K_n}(P_n\times_G M).\end{equation}
This is the simplicial version of the (equivariant) Cartan map. 
If $M$ is a point, this is known as the simplicial Chern-Weil map. 

\subsection{Equivariant Bott-Shulman forms}
Let $G$ be a compact Lie group, and consider the simplicial group 
$E_\bullet G$. The diagonal action 
$ g.(g_0,\ldots,g_n)=(g_0g^{-1},\ldots,g_ng^{-1})$
of $G$ on $E_nG$ makes $E_\bullet G\to B_\bullet G=E_\bullet G/G$ 
into a simplicial principal 
$G$-bundle. 
There is a distinguished 
simplicial connection on $E_\bullet G\to B_\bullet G$, given by 
$$ \sig_n=\sum_{i=0}^n t_i\pr_i^*\theta^L\in 
\Om^1(\Del^n\times E_nG)\otimes\g,$$
where $\pr_i:\,E_nG\to G$ is projection to the $i$th factor and
$\theta^L\in\Om^1(G)\otimes\g$ are the left-invariant Maurer-Cartan
forms. This connection is invariant under the left-action of
$E_\bullet G$ on itself. It hence defines, for any $G$-manifold $M$
and any homomorphism $K_\bullet\to E_\bullet G$ of simplicial groups,  
a $K_\bullet$-equivariant Cartan map \eqref{eq:car1}.  To make this map more
explicit, use the projection
$$G^{n+1}\times M\to G^n\times M,\ \ (g_0,\ldots,g_n,x)\mapsto 
(g_0g_1^{-1},\ldots,g_{n-1}g_n^{-1},g_n.x)$$
to identify $E_nG\times_G M\cong G^n\times M$. Under this 
identification the face maps are 
\beq \partial_i(h_1,\ldots,h_n,x)&=&\left\{\begin{array}{ll}
(h_2,\ldots,h_n,x)&\mbox{for $i=0$},\\
(h_1,\ldots,h_ih_{i+1},\ldots,h_n,x)&\mbox{for $0<i<n$},\\
(h_1,\ldots,h_{n-1},h_n.x)&\mbox{for $i=n$}, 
\end{array}\right.
\eeq
and the action of $E_nG=G^{n+1}$ on $G^n\times M$ reads
$$ (g_0,\ldots,g_n).(h_1,\ldots,h_n,x)=
(g_0 h_1 g_1^{-1},\ g_1 h_2 g_2^{-1},\ldots,g_{n-1}h_n g_n^{-1},g_n.x).$$
We have thus constructed a chain map 
\begin{equation}\label{eq:bosh} 
\phi_{K_\bullet}=\bigoplus_{n\ge 0} \phi^{(n)}_{K_n}:\ \ \Om^\bullet_G(M)\to \bigoplus_{n\ge 0}
\Om^{\bullet-n}_{K_n}(G^n\times M)
\end{equation}
which we will call the {\em equivariant Bott-Shulman map}. As observed
by Bott-Shulman (in a slightly less general setting), the maps
$\phi^{(n)}_{K_n}$ vanishes on $\Om_G^d(M)$ for all $n>d/2$.

To compute the Bott-Shulman map, it is useful to note that the bundle
$E_nG\times M\to E_nG\times_G M$ is trivial: The submanifold defined by 
$g_0=e$ is a trivializing section. In terms of the identification 
$E_nG\times_G M=G^n\times M$, this section is the map 
$$ \iota_n:\,G^n\times M\to G^{n+1}\times M,\ 
(h_1,\ldots,h_n;y)\mapsto (g_0,\ldots,g_n;x)$$
where $g_j=(h_1\cdots h_j)^{-1}$ and $y=(h_1\cdots h_n)^{-1}.x$. 

%

\subsection{The case $M=\pt$}\label{subsec:pt}
The setting originally considered by Bott \cite{bo:le} 
and Shulman \cite{shu:th} corresponds to
the case where $M$ is a point and $K_\bullet$ is trivial. The
equivariant extension to $K_\bullet=G$, viewed as the diagonal
subgroup of $E_\bullet G$, is discussed by Jeffrey in
\cite{je:gr}. In this case, \eqref{eq:bosh} becomes a map
$$ \phi_{G}=\bigoplus_{n\ge 0} \phi^{(n)}_{G}:\ \ 
\on{Pol}^{\bullet}(\g)^G\to \bigoplus_{n\ge 0}
\Om^{2\bullet-n}_{G}(G^n)$$
where $G$ acts on $B_nG=G^n$ by conjugation. To write down concrete 
formulas, note 
that the generating vector field for the diagonal action of $G$ on
$E_nG=G^{n+1}$, given by the left action on each factor, is the 
sum of $n+1$ copies of the right-invariant vector field 
$-\xi^R$. Hence, the equivariant curvature of the connection $\sig_n$ reads 
$$ F_G(\xi)=\sum_{i=0}^n \d t_i\, g_i^*\theta^L-\hh \sum_{i=0}^n t_i\,g_i^*[\theta^L,\theta^L]
-\hh \sum_{i,j=0}^n t_i t_j [g_i^*\theta^L,g_j^*\theta^L]
+\sum_{i=0}^n t_i\,\Ad_{g_i^{-1}}(\xi).$$
Then $\phi^{(n)}(p)$ is given explicitly as the integral 
$$\phi^{(n)}(p)=\int_{\Del^n} p(\iota_n^* F_G(\xi)).$$ 
For $p\in \on{Pol}^d(g)^G$, introduce special notation 
\beq \eta^p_G&=&\phi^{(1)}_G({p})\in\Om^{2d-1}_G(G),\\
\varphi^p&=&\phi^{(2)}_G({p}) \in \Om^{2d-2}_G(G^2)\eeq
(obviously $\phi^{(0)}_G({p})=0$ if $d>0$). If $p(\xi)=\hh
||\xi||^2$ we will drop the superscript $p$. It is not hard to see that 
the general formula for $\phi^{(n)}(p)$ specializes, for $n=1$, to 
the formula \eqref{eq:jef} for $\eta^p_G$. Also, taking $n=2$ 
and  $p(\xi)=\hh ||\xi||^2$ one finds
that $\varphi$ is an ordinary 2-form on $G\times G$:
$$\varphi=\hh \pr_1^*\theta^L\cdot \pr_2^*\theta^R.$$
(The forms $\phi^{(n)}_G({p})$ for $n>2$ vanish for $p(\xi)=\hh
||\xi||^2$.)  

The fact that $\phi_G({p})$ is closed under the total differential for
the double complex $(\Om^{k}_G(G^l),\d_G,\delta)$ gives equations,
\begin{eqnarray} 
 \d_G\eta^p_G&=&0\nonumber,\\ 
\d_G\varphi^p&=&\on{Mult}^*\eta^p_G-\pr_1^*\eta^p_G-\pr_2^*\eta^p_G
\label{eq:prim}
\end{eqnarray}
where $\on{Mult}:\,G\times G\to G$ is group multiplication. 
Furthermore, the form
\beq
\lefteqn{(\pr_1\times \pr_2)^*\varphi^p+ (\on{Mult}\circ (\pr_1\times
\pr_2)\times \pr_3)^*\varphi^p}\\&&
-(\pr_2\times\pr_3)^*\varphi^p-(\pr_1\times (\on{Mult}\circ
(\pr_2\times \pr_3))^*\varphi^p
\eeq
on $G\times G\times G$ is $\d_G$-exact. Another useful property of the forms $\eta^p_G$ is that they change sign under the inversion map 
$\on{Inv}:\,G\to G,\ g\mapsto g^{-1}$ 
\begin{equation}\label{eq:inversion}
 \on{Inv}^*\eta_G^p=-\eta_G^p.
\end{equation}
This follows from Remark 2.1(b), by pulling the identity
\eqref{eq:prim} back under the map $g\mapsto (g,g^{-1})$, and using
that the pull-back of $\eta_G^p$ to the group unit vanishes.

While the formulas for the forms $\varphi^p$ are rather complicated, 
one has simple expressions if $G=T$ is a torus. 

\begin{lemma}\label{lem:varphiab}
If $G=T$ is a torus, the form $\varphi^p(\xi)\in \Om_T(T^2)$ is given by 
$$ 
\varphi^p(\xi)=\hh\sum_{rs} p''_{rs}(\xi)\pr_1^*\theta_T^r \pr_2^*\theta_T^s.
$$
\end{lemma}
\begin{proof}
In the Abelian case, the formula for the equivariant curvature of 
$\sig_2$ simplifies dramatically:
$$ F_T(\xi)=\sum_{i=0}^2 \d t_i\, g_i^*\theta_T+\xi.$$
Hence, the pull-back under $\iota_2(h_1,h_2)=(e,h_1^{-1},(h_1h_2)^{-1})$ 
reads, 
$$ \iota_2^*F_T(\xi)=\xi-\d t_1\ h_1^*\theta_T -\d t_2 (h_1^*\theta_T+h_2^*\theta_T).$$
Hence, 
\beq p(\iota_2^*F_T(\xi))
&=&p(\xi -\d t_1\ h_1^*\theta_T -\d t_2 (h_1^*\theta_T+h_2^*\theta_T))\\
&=&-\d t_1 \d t_2 \sum_{rs}\,p''_{rs}(\xi)\, h_1^*\theta_T^r\, h_2^*\theta_T^s
+\ldots
\eeq 
where we have only written the coefficient of $\d t_1 \d t_2$. The 
Lemma follows since $\int_{\Del^2}\d t_1 \d t_2=-\hh$, with our conventions. 
\end{proof}

\subsection{The case $M=G$}\label{subsec:G}
Consider next the case $M=G$ with $G$ acting by conjugation, and with
$K_\bullet=E_{\bullet}G=G^{\bullet+1}$. We obtain a chain map
$$ \psi=\bigoplus_{n=0}^\infty \psi^{(n)}:\ \ 
\Om^{\bullet}_G(G)\to \bigoplus_{n=0}\Om^{\bullet-n}_{G^{n+1}}(G^{n+1})$$
As pointed out above, $\psi^{(0)}$ is just the identity map. 
Consider the degree $n=1$ component. 
The action of $G^2$ is given by $(g_0,g_1).(h,k)=(g_0 h
g_1^{-1},g_1 k g_1^{-1}) $, and the face maps are 
$\partial_0(h,k)=k$ and $\partial_1(h,k)=\Ad_h(k)$. 
The form 
$$ \lambda^p:=\psi^{(1)}(\eta_G^p)\in\Om^{2d-2}_{G^2}(G^2)$$ 
has the property
\begin{equation}\label{eq:lambda1}
 \textstyle{\d_{G^2}\lambda^p(\xi_0,\xi_1)=\partial_0^*\eta_G^p(\xi_1)-
\partial_1^*\eta_G^p(\xi_0)}.
\end{equation}
Again we drop the superscript for $p(\xi)=\hh||\xi||^2$, 
and also the subscript $G^2$ since  $\lambda^p$ does 
not depend on the equivariant parameter in this case. 
Explicit calculation gives:
\begin{equation}\label{eq:lambda}
\lambda=-\hh h^*\theta^L\cdot \Ad_k(h^*\theta^L)+
h^*\theta^L\cdot k^*\f{\theta^L+\theta^R}{2}
\end{equation}
(where we view $h,k$ as maps $G^2\to G$). For more general 
$p$, we only state the result if $G=T$ is a torus: 

\begin{lemma}
If $G=T$ is a torus, 
$$ \lambda^p(\xi_0,\xi_1)=\sum_{rs} A_{rs}(\xi_0,\xi_1)
 h^*\theta_T^r\, k^*\theta_T^s$$
where $A_{rs}(\xi_0,\xi_1)=\int_0^1 \,p''_{rs}(t\xi_0+(1-t)\xi_1)\ \d t$. 
\end{lemma}
\begin{proof}
We denote points in $E_nT\times T=T^{n+1}\times T$ by
$(g_0,\ldots,g_n,k)$, and points in $E_nT\times_T G=T^n\times T$ by
$(h_1,\ldots,h_n,k)$.  The $E_nT=T^{n+1}$-equivariant curvature of $\sig_n$ is
given by
$$ F^\sig_{T^{n+1}}(\xi_0,\ldots,\xi_n)=\sum_{i=0}^n t_i \xi_i+ \sum_{i=0}^n
\d t_i\ g_i^*\theta_T
$$
Recall that $\eta^p_T(\xi)=-p'(\xi)\cdot\theta_T$. 
To compute $\lambda^{(p)}=\psi^{(1)}(\eta^p_T)$, we have to consider the form
in $\Om_{T^2}(E_1T\times T)$, 
$$ \eta^p(F^\sig_{T^2}(\xi_0,\xi_1))=-
k^*\theta_T\cdot p'(\sum_{i=0}^1 t_i \xi_i- \sum_{i=0}^1 \d t_i\ g_i^*\theta_T)
$$
Pulling back under the map $\iota_1(h,k)=(e,h^{-1},k)$, and setting $t_1=t,\,
t_0=(1-t)$ we obtain
\beq \iota_1^*\eta^p(F^\sig_{T^2}(\xi_0,\xi_1))&=&-
k^*\theta_T\cdot p'((1-t)\xi_0+t\xi_1+\d t h^*\theta_T)\\
&=& \d t \sum_{rs} p''_{rs}((1-t)\xi_0+t\xi_1)
k^*\theta_T^r h^*\theta_T^s+\ldots
\eeq
where we have only written the coefficient of $\d t$. Integrating over 
$\Del^1$, the Lemma follows. 
\end{proof}

\end{appendix}


\def\polhk#1{\setbox0=\hbox{#1}{\ooalign{\hidewidth
  \lower1.5ex\hbox{`}\hidewidth\crcr\unhbox0}}} \def\cprime{$'$}
  \def\cprime{$'$} \def\polhk#1{\setbox0=\hbox{#1}{\ooalign{\hidewidth
  \lower1.5ex\hbox{`}\hidewidth\crcr\unhbox0}}} \def\cprime{$'$}
\providecommand{\bysame}{\leavevmode\hbox to3em{\hrulefill}\thinspace}
\providecommand{\MR}{\relax\ifhmode\unskip\space\fi MR }
\providecommand{\MRhref}[2]{%
  \href{http://www.ams.org/mathscinet-getitem?mr=#1}{#2}
}
\providecommand{\href}[2]{#2}

\end{document}